\documentclass[11pt, oneside]{amsart}
\usepackage{amsfonts, amscd, amsmath, amssymb}
\usepackage[dvips]{graphicx}

\begin{document}

\newcommand\globcnt{subsection}
\newcommand\globcntthe{\thesubsubsection}
\newcommand\globsect{\subsubsection*}

\newcommand\MainTheorem{Main~Theorem}

\theoremstyle{plain}
\newtheorem*{MainTh}{\MainTheorem}
\newtheorem*{condition}{Condition}
\newtheorem{thm}[\globcnt]{Theorem}
\newtheorem{lem}[\globcnt]{Lemma}
\newtheorem{claim}[\globcnt]{Claim}
\newtheorem{prop}[\globcnt]{Proposition}
\newtheorem{cor}[\globcnt]{Corollary}
\newtheorem{add}[\globcnt]{Addendum}
\theoremstyle{definition}
\newtheorem{defn}[\globcnt]{Definition}
\newtheorem{exmp}[\globcnt]{Example}
\newtheorem{rem}[\globcnt]{Remark}

\makeatletter
\@addtoreset{subsubsection}{section}
\@addtoreset{subsubsection}{subsection}
\@addtoreset{equation}{subsection}
\@addtoreset{figure}{subsection}
\makeatother

\renewcommand{\theequation}{\thesection.\arabic{equation}}
\renewcommand{\thefigure}{\thesection.\arabic{figure}}

\newcommand\chEq{%
}

\newcommand\chFig{%
}

\newcommand\CCC{{\mathbb C}}
\newcommand\NNN{{\mathbb N}}
\newcommand\PPP{{\mathbb P}}
\newcommand\RRR{{\mathbb R}}
\newcommand\SSS{{\mathbb S}}
\newcommand\ZZZ{{\mathbb Z}}

\newcommand\AAA{{\mathcal A}}
\newcommand\BB{{\mathcal B}}
\newcommand\CC{{\mathcal C}}
\newcommand\DD{{\mathcal D}}
\newcommand\EE{{\mathcal E}}
\newcommand\FF{{\mathcal F}}
\newcommand\GG{{\mathcal G}}
\newcommand\II{{\mathcal I}}
\newcommand\JJ{{\mathcal J}}
\newcommand\KK{{\mathcal K}}
\newcommand\MM{{\mathcal M}}
\newcommand\NN{{\mathcal N}}
\newcommand\OO{{\mathcal O}}
\newcommand\PP{{\mathcal P}}
\newcommand\RR{{\mathcal R}}
\newcommand\TT{{\mathcal T}}
\newcommand\UU{{\mathcal U}}
\newcommand\VV{{\mathcal V}}
\newcommand\WW{{\mathcal W}}
\newcommand\XX{{\mathcal X}}
\newcommand\YY{{\mathcal Y}}
\newcommand\ZZ{{\mathcal Z}}

\newcommand\Int{\mathrm{Int}}
\newcommand\id{\mathrm{id}}
\newcommand\grad{\mathrm{grad}}
\newcommand\supp{\mathrm{supp\ }}
\newcommand\IM{\mathrm{Im}}
\newcommand\eps{\varepsilon}

\newcommand\surf{M} 
\newcommand\manif{\surf}
\newcommand\bdsurf{\surf}
\newcommand\ubdsurf{\widehat\surf}
\newcommand\shmanif{\widehat{\manif}}

\newcommand\Diff{\mathcal{D}}            
\newcommand\DiffId{\Diff_{\id}}
\newcommand\DiffM{\Diff(\manif)}

\newcommand\Fld{F}

\newcommand\pnt{z}
\newcommand\hommap{*}

\newcommand\onemanif{P}
\newcommand\Circle{S^1}
\newcommand\Rline{\RRR^1}
\newcommand\MrsMP{\FF(\surf,\onemanif)}
\newcommand\MrsMR{\FF(\surf,\Rline)}
\newcommand\MrsMS{\FF(\surf,\Circle)}

\newcommand\smsp[1]{C^{\infty}\ifx#1\empty\else(#1)\fi}
\newcommand\smMP{\smsp{\surf,\onemanif}}
\newcommand\smMR{\smsp{\surf,\Rline}}
\newcommand\smMS{\smsp{\surf,\Circle}}

\newcommand\func{f}
\newcommand\gfunc{g}
\newcommand\hdif{h}
\newcommand\pdif{p}
\newcommand\thdif{\widetilde{\hdif}}
\newcommand\adif{p}
\newcommand\bdif{q}

\newcommand\Hisot{H}
\newcommand\Asgm{\Phi}
\newcommand\Bsgm{\Psi}
\newcommand\tAsgm{\widetilde{\Asgm}}

\newcommand\barfunc{\bar\func}

\newcommand\Fsgm{F}
\newcommand\tFsgm{\Fsgm}
\newcommand\fsgm{\func}
\newcommand\tfsgm{\tfunc}

\newcommand\Gsgm{G}

\newcommand\crpnt[2]{c_{#1}(#2)}
\newcommand\crpntf[1]{\crpnt{#1}{\func}}
\newcommand\crpntg[1]{\crpnt{#1}{\gfunc}}
\newcommand\crpt[1]{c_{#1}}

\newcommand\difr[2]{d#1(#2)}
\newcommand\tansp[2]{T#1_{#2}}

\newcommand\bcomp{V}
\newcommand\sgn[1]{\eps_{#1}}
\newcommand\bsgn[1]{\eps_{#1}}

\newcommand\Crtype[1]{K(#1)}
\newcommand\CrtypeEmpty{K}

\newcommand\tmpcurve{\gamma}
\newcommand\ttmpcurve{\widetilde{\tmpcurve}}
\newcommand\sgmhom{$\Sigma$}
\newcommand\sgmdef[1]{\, \stackrel{#1}{\sim} \, }
\newcommand\sgmh{\stackrel{\Sigma}{\sim}}

\newcommand\Dtwist{t}
\newcommand\tDtwist{\widetilde{\Dtwist}}
\newcommand\Dtwcurv{\Dtwist_{\tmpcurve}}
\newcommand\tDtwcurv{\Dtwist_{\tcurve}}
\newcommand\Dtwc[1]{\Dtwist_{\tmpcurve_{#1}}}
\newcommand\Dtw[1]{\Dtwist_{#1}}

\newcommand\KrReeb{\Gamma}
\newcommand\ReebGr[1]{\KrReeb_{#1}}
\newcommand\Reebf{\ReebGr{\func}}
\newcommand\Reebg{\ReebGr{\gfunc}}

\newcommand\pReebf{\func^{*}}
\newcommand\funcR{\func_{\KrReeb}}
\newcommand\pReebg{\gfunc^{*}}
\newcommand\gfuncR{\gfunc_{\KrReeb}}

\newcommand\edgeR{e}

\newcommand\bNbh{N}

\newcommand\Rdif{\alpha}
\newcommand\Cdif{\phi}
\newcommand\Cisot{\Phi}

\newcommand\generS{\eta}
\newcommand\imgenerS{\xi}

\newcommand\intform[2]{
\omega
\ifx#1\cdot
  \ifx#2\cdot
  \else\fi
\else
  \left(#1,#2\right)
\fi
}

\newcommand\acurve{\alpha}
\newcommand\bcurve{\beta}
\newcommand\ccurve{\gamma}
\newcommand\dcurve{\delta}
\newcommand\ecurve{\epsilon}
\newcommand\scurve{\sigma}

\newcommand\aacurve{\hat\acurve}

\newcommand\tbcurve{\widetilde{\beta}}
\newcommand\tdcurve{\widetilde{\delta}}

\newcommand\acoef{a}
\newcommand\bcoef{b}

\newcommand\germ{g}
\newcommand\hdifhom{\hdif_{\hommap}}
\newcommand\Sp[2]{Sp_{#1}(#2)}

\newcommand\isointers{\phi}

\newcommand\Idmatr{I}
\newcommand\Eijmatr{e_{ij}}

\newcommand\Ematr[2]{e_{#1#2}}

\newcommand\aastddif[2]{\mu_{#1#2}}
\newcommand\abstddif[2]{\nu_{#1#2}}
\newcommand\bbstddif[2]{\eta_{#1#2}}

\newcommand\comp{X}
\newcommand\collar{N}

\newcommand\cov{\sigma}
\newcommand\indk{k}

\newcommand\sVar{s}

\newcommand\posbnd{\partial_{+}\surf}
\newcommand\negbnd{\partial_{-}\surf}

\newcommand\posbndf{\partial_{+}\surf}
\newcommand\negbndf{\partial_{-}\surf}

\newcommand\tposbndf{\partial_{+}\tsurf}
\newcommand\tnegbndf{\partial_{-}\tsurf}

\newcommand\tposbndfX{\partial_{+}\comp}
\newcommand\tnegbndfX{\partial_{-}\comp}

\newcommand\cval{c}
\newcommand\aval{a}

\newcommand\tCircle{\widetilde{S}}
\newcommand\tfunc{\widetilde{\func}}
\newcommand\tgfunc{\widetilde{\gfunc}}
\newcommand\hfunc{\widehat{\func}}
\newcommand\hgfunc{\widehat{\gfunc}}

\newcommand\tsurf{\widetilde{\surf}}
\newcommand\degi{d}
\newcommand\psurf{p}
\newcommand\pcirc{q}
\newcommand\tcurve{\widetilde{\tmpcurve}}
\newcommand\btsurf[1]{B_{#1}}

\newcommand\tCrv{K}
\newcommand\grpAA{\tCrv_{0}}
\newcommand\grpBB{\tCrv^{1}}
\newcommand\grpAB{\tCrv_{0}^{1}}

\newcommand\tComp{Q}
\newcommand\compAA{\tComp_{0}}
\newcommand\compBB{\tComp^{1}}
\newcommand\compAB{\tComp_{0}^{1}}

\newcommand\lowarcs{G}

\newcommand\Nbh{V}
\newcommand\aNbh{W}

\newcommand\divcurves{Z}
\newcommand\mincurves{\widetilde{\divcurves}}
\newcommand\arc{l}

\newcommand\Invim{Z}
\newcommand\surfp{V}

\newcommand\onehalf{\frac{1}{2}}
\newcommand\Onehalf{1/2}

\newcommand\oneforth{\frac{1}{4}}
\newcommand\threeforth{\frac{3}{4}}

\newcommand\bndpos[1]{b_{+}(#1)}
\newcommand\bndneg[1]{b_{-}(#1)}

\newcommand\bndposf{b_{+}}
\newcommand\bndnegf{b_{-}}

\newcommand\tbndposf{\bndpos{\tfunc}}
\newcommand\tbndnegf{\bndneg{\tfunc}}

\newcommand\pFld{\Phi}
\newcommand\gFld{\Psi}

\newcommand\tFld{\widetilde{\Fld}}
\newcommand\tpFld{\widetilde{\pFld}}
\newcommand\tgFld{\widetilde{\gFld}}

\newcommand\canonical{canonical}

\newcommand\frval[2]{(#1,#2)}
\newcommand\concomp[1]{\#[#1]}

\newcommand\AdmDif[1]{\AAA(#1)}
\newcommand\PNDif[1]{NP(#1)}
\newcommand\AdmOrb[1]{\OO(#1)}
\newcommand\tcomp{\widetilde{\comp}}
\newcommand\PathComp[1]{C(#1)}
\newcommand\fAdmDif{\AdmDif{\func}}
\newcommand\gAdmDif{\AdmDif{\gfunc}}

\newcommand\mybrack[1]{\left\| #1 \right\|}

\newcommand\Torelly[1]{\mathcal{T}(#1)}
\newcommand\sTrl{\Torelly{\surf}}

\newcommand\cnmap[1]{\widehat{#1}}
\newcommand\cnfunc{\cnmap{\func}}
\newcommand\cngfunc{\cnmap{\gfunc}}

\newcommand\MapClGr{\MM}
\newcommand\HomGr[2]{\MapClGr_{#2}(#1)}
\newcommand\PHomGr[2]{\PP\!\MapClGr_{#2}(#1)}
\newcommand\HM{\HomGr{\surf}{}}
\newcommand\HMf{\HomGr{\surf}{\func}}
\newcommand\HMbndpl{\MapClGr^{+}_{\partial}(\surf)}
\newcommand\TrlGrM{\mathcal{T}(\surf)}
\newcommand\TrlGrMu{\mathcal{T}(\ubdsurf)}
\newcommand\HMn{\HomGr{\ubdsurf\setminus\{x_1,\ldots,x_n\}}{}}

\newcommand\HMnubd{\HomGr{\ubdsurf}{n}}
\newcommand\HMbd{\HomGr{\bdsurf}{}}
\newcommand\HMubd{\HomGr{\ubdsurf}{}}
\newcommand\PHMbd{\PP\!\HomGr{\bdsurf}{}}

\newcommand\CurvesSet{S}

\newcommand\sgmdif{\sigma}
\newcommand\rhodif{\rho}
\newcommand\bndif{b}

\newcommand\DifOr{O}

\newcommand\ConfSp[2]{F_{#1}(#2)}
\newcommand\clsurf{\surf'}

\newcommand\Dnsurf{\Diff(\ubdsurf)_{n}} 
\newcommand\DnsurfId{\Diff_{\id}(\ubdsurf)_{n}} 

\newcommand\Dfsurf{\Diff\surf}
\newcommand\Dudbsurf{\Diff\ubdsurf}

\newcommand\embmap{*}
\newcommand\inclhom{\stackrel{\embmap}{\hookrightarrow}}

\newcommand\cldif{c}
\newcommand\brdif{b}
\newcommand\trldif{t}
\newcommand\orpermdif{p}

\newcommand\SpZg{\Sp{2\germ}{\ZZZ}}
\newcommand\modZg{\ZZZ^{2\germ}}
\newcommand\Ann{A}
\newcommand\TransAnn{T}
\newcommand\Stab{\mathrm{St}}
\newcommand\Hab{H}

\newcommand\CzMS{C^{0}(\surf,\Circle)}
\newcommand\HclMS{[\surf,\Circle]}
\newcommand\MHclMS{\MM[\surf,\Circle]}
\newcommand\bhh{\eta}
\newcommand\lset{L}

\newcommand\acrc{l}
\newcommand\narc{k}
\newcommand\arccompE{L_{\varnothing}}
\newcommand\arccompAA{L_0}
\newcommand\arccompAB{L_0^1}
\newcommand\arccompBB{L^1}
\newcommand\indnum{n}
\newcommand\Arcs{L}
\newcommand\sumdeg{d}
\newcommand\gunion{\Gamma}
\newcommand\genus{g}

\newcommand\KR{$K\!R$}
\newcommand\rankH{r}
\newcommand\obr{q}
\newcommand\mystar{*}

\newcommand\LS{$LS$}
\newcommand\Syst{\Gamma}
\newcommand\mcut[1]{\surf_{#1}}
\newcommand\fcut[1]{\func_{#1}}
\newcommand\gcut[1]{\gfunc_{#1}}
\newcommand\prcut[1]{p_{#1}}
\newcommand\compall[1]{B_{#1}}
\newcommand\compneg[1]{\compall{#1}^{-}}
\newcommand\comppos[1]{\compall{#1}^{+}}
\newcommand\mcomp{Q}
\newcommand\mcomppos{\mcomp^{+}}
\newcommand\mcompneg{\mcomp^{-}}
\newcommand\mcompboth{\mcomp^{\pm}}
\newcommand\Curv{\Lambda}

\newcommand\intrsc[2]{\langle#1,#2\rangle}
\newcommand\Ta{\TT_{\alpha}}
\newcommand\Tam{\Ta(\surf)}
\newcommand\Tama{\UU}
\newcommand\Torel{\TT(\surf)}
\newcommand\hacurve{\acurve'}
\newcommand\aafunc{\func}
\newcommand\typeC{(c)}
\newcommand\gdif{g}
\newcommand\nacurve{\acurve_1}
\newcommand\larc{l}

\newcommand\ydif{y}
\newcommand\tydif{\widetilde{\ydif}}
\newcommand\Kleinbtl{K}
\newcommand\Torhol{T}
\newcommand\Mband{B}
\newcommand\Mobius{M\"obius}

\newcommand\SCC{SCC}
\newcommand\SCCs{SCCs}
\newcommand\bsldif{\nu}
\newcommand\obsldif{\omega}

\newcommand\vval{v}
\newcommand\wval{w}

{\centering
{\Large\bf
Path-components of Morse mappings spaces of surfaces
}\\

\medskip

{\sc Sergey Maksymenko\footnote{The author is partially supported by the grant of Government Fond of Fundamental Researches no. 1.7/132 }}\\

\medskip

Topology Dept.,
Institute of Mathematics of Ukrainian NAS,
Te\-re\-shchen\-kiv\-ska str. 3,
01601 Kyiv,
Ukraine,
e-mail:~\texttt{maks@imath.kiev.ua}

\bigskip
\smallskip

\begin{minipage}[c]{0.9\textwidth}
\small
\centerline{\bf Abstract}
Let $\surf$ be a connected compact surface, $\onemanif$ be either $\Rline$ or $\Circle$, and
$\MrsMP$ be the space of Morse mappings $\surf\to\onemanif$ with compact-open topology.
The classification of path-components of $\MrsMP$
was independently obtained by S.~V.~Matveev and
V.~V.~Sharko for the case $\onemanif=\Rline$, and by the author for orientable surfaces and $\onemanif=\Circle$.
In this paper we give a new independent and unified proof of this classification for all compact surfaces in the case $\onemanif=\RRR$,
and for orientable surfaces in the case $\onemanif=\Circle$.
We also extend the initial author's proof to non-orientable surfaces.
\end{minipage}

}

\bigskip
\smallskip


\section{Introduction}
Let $\surf$ be a smooth ($C^{\infty}$) connected compact surface with boundary $\partial\surf$ (possibly empty) and $\onemanif$ be a one-dimensional manifold, i.e. either the real line $\Rline$ or the circle $\Circle$.
Consider the subspace $\MrsMP$ of $\smMP$ consisting of Morse mappings $\surf\to\onemanif$.
It is well-known that $\MrsMP$ is an everywhere dense open subset of $\smMP$ in the compact-open topology of $\smMP$.
The homotopy type of this space is of great importance in differential topology and dynamical systems, see e.g.~\cite{Hirsh:DT, Igusa, ThurstonHatcher, HirshHirsh, Kudr, Sharko, Maks:PhD, SaekiIkegami}.

Recently, S.~V.~Matveev and V.~V.~Sharko~\cite{Sharko} have obtained a full description of path-components of the space $\MrsMR$. 
Matveev's proof is included and generalized in the paper~\cite{Kudr} of E.~Kud\-ryavt\-seva to numerated Morse functions. 
Their proofs were independent and based on different ideas.
The classification of path-components of $\MrsMS$ for orientable surfaces was given in the author's Ph.D., see~\cite{Maks:PhD}. 

These results (which we will refer to as \MainTheorem) can be summarized as follows:
two Morse mappings $\func,\gfunc:\manif\to\onemanif$ belong to same path-component of $\MrsMP$ if and only if they are homotopic as continuous maps and have the same number of critical points at each index and the same sets of {\em positive\/} and {\em negative\/} boundary components (in the sense described below.)

In this paper we give a unified and independent proof of this theorem for all compact surfaces in the case $\onemanif=\RRR$.
The case of Morse mappings $\surf\to\Circle$ requires information on the subgroup of the mapping class group of $\surf$ preserving a given element in the cohomology group $H^1(\surf,\ZZZ)$.
We also find the generators of this group for orientable surfaces and 
extend the presented method to Morse mappings from orientable surfaces into $\Circle$.

In fact, the proof given in~\cite{Maks:PhD} for this case almost literally extends to non-orientable surfaces as well. Since~\cite{Maks:PhD} was never published in English, we give this proof for all surfaces in Appendix.
Thus the \MainTheorem\ is proved here for all cases of $\surf$ and $\onemanif$.

Our approach has a relation to the paper~\cite{ThurstonHatcher} of A.~Hatcher and W.~Thurs\-ton, who used deformations of Morse functions to construct a representation for the mapping class group of a surface.
In contrast to this approach, we exploit generators of this group to find a deformation between Morse mappings in $\MrsMP$.
The principal observation is that ``elementary diffeomorphisms'' like Dehn twists, boundary and crosscap slides generating mapping class groups of surfaces preserve certain Morse functions.


\section{Preliminaries}\label{sect:Constructions}
Let $\manif$ be a compact surface.
A surface obtained by shrinking every connected component of $\manif$ to a point will be denoted by $\shmanif$.
Thus $\shmanif$ is closed and is homeomorphic with a connected sum of the form either
$S^2 \mathop\#\limits_{i=1}^{\genus} T^2$ (orientable case, $\genus\geq 0$) or
$\mathop\#\limits_{i=1}^{\genus} \PPP^2$ (non-orientable case, $\genus\geq 1$).
In each of the cases the number $\genus$ is called the germ of $\manif$.
All homology and cohomology groups will be taken with integer coefficients.
The term {\em simple closed curve} will be abbreviated to \SCC.
The circle $\Circle$ will be regarded as the subset $\{z\in\CCC: |z|=1 \}$ of the complex plane $\CCC$.
For a topological space $X$ let $\concomp{X}$ denote the number of its connected components.

\subsection{Morse mappings}
Let us fix, once and for all, an orientation of $\onemanif$.
Consider a smooth mapping $\func:\surf\to\onemanif$.
A point $\pnt\in\surf$ is {\em critical} for $\func$ if $\difr{\func}{\pnt}=0$.
A critical point $\pnt$ of $\func$ is {\em non-degenerate} if the Hessian
of $\func$ at $\pnt$ is non-degenerate.
Suppose that $\pnt$ is a non-degenerate critical point of $\func$.
Then by Morse lemma there are embeddings $\psurf:\RRR^2\to\surf$ and $\pcirc:\Rline\to\onemanif$
onto open neighborhoods of $\pnt$ and $\func(\pnt)$ respectively
such that $\psurf(0,0)=\pnt$, $\pcirc(0)=\func(\pnt)$,
$\pcirc$ preserves orientation, and $\pcirc^{-1}\circ\func\circ\psurf(x,y)=\pm x^2 \pm y^2$.
The number of minuses in this representation
does not depend on a particular choice of such embeddings and
is called the {\em index\/} of a critical point $\pnt$.

A $C^{\infty}$-mapping $\func:\surf\to\onemanif$ is {\em Morse\/} if the following conditions hold:
\begin{enumerate}
\item
all critical points of $\func$ are non-degenerate and belong to the interior of $\surf$;
\item
$\func$ is constant at each boundary component of $\surf$
while its values on different components may differ each from other.
\end{enumerate}
The subspace of $\smMP$ consisting of Morse mappings will be denoted by $\MrsMP$.
We endow $\smMP$ with the compact-open topology. Then this topology induces
some topology on $\MrsMP$.

\subsection{\sgmhom-homotopies.}
Let $\func,\gfunc\in\MrsMP$ be two Morse mappings
and $\phi:[0,1]\to\smMP$ be a path between them in the space of Morse mappings, thus $\phi$ is continuous, $\phi(0)=\func$, $\phi(1)=\gfunc$ and $\phi(t)$ is Morse for all $t\in[0,1]$.
Then $\phi$ yields a continuous mapping (homotopy)
$F:\surf\times I\to \onemanif$ such that
$F_{0} = \func$, $F_{1} = \gfunc$, and $F_{t}$ is Morse for all $t\in I$.
In particular, $F$ is $C^{\infty}$ in $x\in\surf$ but may be just continuous in $t\in[0,1]$.
Conversely, every such mapping $F$ gives rise to a path between $\func$ and $\gfunc$ in $\MrsMP$.

We will call the mapping $F$ a {\em \sgmhom-homotopy} or {\em \sgmhom-deformation\/}
between $\func$ and $\gfunc$
and write $\func \sgmdef{F_t} \gfunc$.
The term $\func \sgmh \gfunc$ will also be used to indicate
that $\func$ and $\gfunc$ are \sgmhom-homotopic.

\begin{rem}
In~\cite{Sharko, Kudr} \sgmhom-homotopies are called {\em isotopies\/} of Morse functions.
We will use another term in order to avoid confusions with isotopies of diffeomorphisms.
\end{rem}

\subsection{Invariants of \sgmhom-homotopies.}
Let $\func\in\MrsMP$.
The objects (i) homotopy class, (ii) number of critical points in each index, and (iii) positive and negative boundary components are invariant under \sgmhom-homotopies of $\func$.


\subsubsection{Homotopy class.}
First suppose that $\onemanif=\Circle$.
Let $\xi\in H^{1}(\Circle)$ be a generator defining the chosen orientation of $\Circle$.

If $\func:\surf\to\Circle$ is a continuous mapping, then the correspondence $\func\mapsto\func^{*}(\xi)\in\ H^{1}(\surf)$ yields a bijection between the set of homotopy classes of mappings $[\surf,\Circle]$ and the cohomology group $H^{1}(\surf)$. 
Since by our definition Morse mappings are constant at the connected components of $\surf$, it follows that the set of homotopy classes of Morse mappings $\surf\to\Circle$ is bijective to the group $H^{1}(\shmanif)$ for the corresponding closed surface $\shmanif$.


Let $\genus$ be a genus of $\manif$.
A simple calculation shows that $H^{1}(\shmanif)$ is isomorphic with $\ZZZ^{\rankH}$, where $\rankH$ is either $2g$ or $g-1$ provided $\manif$ is orientable or not.
Let us fix a basis for $H^{1}(\shmanif)$. 
Then the homotopy class of $\func$ is an integer vector 
$$ (\obr_1,\ldots,\obr_{\rankH}) = \func(\xi) \in H^{1}(\shmanif)=\ZZZ^{\rankH}.$$
For $\onemanif=\RRR$ we will assume that $(\obr_1,\ldots,\obr_{\rankH})=(0,\ldots,0)$.

\subsubsection{Number of critical points in each index.}
Denote by $\crpntf{i}=\crpt{i}$, $(i=0,1,2)$ the number of critical points of $\func$ of index $i$.
Then by Morse equalities we have
\chEq\begin{equation}\label{equ:Morse_equ}
\crpntf{0} + \crpntf{1} - \crpntf{2} = \chi(\manif).
\end{equation}

\subsubsection{Positive and negative components of $\partial\manif$.}
Let $\bcomp$ be a component of $\partial\surf$,
$\pnt\in\bcomp$ and $\xi\in\tansp{\surf}{\pnt}$ be a tangent vector at $\pnt$
directed {\em outward\/} $\surf$.
Denote by $\bsgn{\func}(\bcomp)$ the {\em sign\/} of the value $\difr{\func}{\pnt}\xi$.
Since $\func$ has no critical points on $\bcomp$,
we see that $\bsgn{\func}(\bcomp)=\pm1$ and does not depend on a particular choice of
a point $\pnt\in\bcomp$ and a vector $\xi\in\tansp{\surf}{\pnt}$ as above.
Thus we get a function $\bsgn{\func}:\pi_0\partial\manif\to\{\pm 1\}$.
We may also think of $\bsgn{\func}$ as an element of $\{\pm 1\}^{b}$, where $b$ is the number of connected components of $\partial\manif$.

We will call $\bcomp$ either {\em $\func$-positive\/} or {\em $\func$-negative\/}
in accordance with $\bsgn{\func}(\bcomp)$.
Let $\posbndf$ (resp. $\negbndf$) be the union of $\func$-positive (resp. $\func$-negative) boundary components of $\partial\surf$, and let $b_{+}$ ($b_{-}$) denote the numbers of these components.

The following collection of numbers
$$
 \Crtype{\func} =
 \left\{
\obr_{1},\ldots,\obr_{\rankH}, \
\crpt{0}, \crpt{1} , \crpt{2}, \
\bsgn{\func}
 \right\}
$$
will be called the {\em critical type\/} of a Morse mapping $\func$.
It can be regarded as a point in
$\ZZZ^{\rankH} \times \NNN_{0}^{3} \times \{\pm1\}^{b}$
belonging to the ``hyperplane'' defined by Eq.~\eqref{equ:Morse_equ}, where $\NNN_{0}=\NNN\cup\{0\}$.
If we choose another orientation of $\onemanif$, then $\crpntf{0}$ exchanges with $\crpntf{2}$, $\crpntf{1}$ remains unchanged, $\bsgn{\func}$ and every $\obr_{i}$ change their signs.

Our aim is to give a new proof of the following theorem:
\begin{MainTh}[Matveev~\cite{Kudr}, Sharko~\cite{Sharko}, Maksymenko~\cite{Maks:PhD}]
Two Morse mappings $\func, \gfunc:\surf\to\onemanif$ belong to the same path-component of $\MrsMP$ if and only if $\Crtype{\func} = \Crtype{\gfunc}$, i.e.\! they are homotopic, have the same numbers of critical points in each index, and the same sets of positive and negative components of $\partial\manif$.
\end{MainTh}
The necessity is obvious therefore we confine ourself by the sufficiency.
Let us briefly review the known proofs of this theorem.
First consider the case $\onemanif=\Rline$.
Let $\func$ and $\gfunc$ be two Morse functions with equal critical types.
In the both proofs~\cite{Kudr,Sharko} the problem was reduced to minimal Morse functions with no critical points of index $0$ and $2$.

Let $\Fld$ be a gradient-like vector field for a minimal Morse
function $\func$. Consider a union of $\func$-negative boundary
components of $\surf$ with trajectories of $\Fld$ that finish at
critical points of $\func$.
This set is called a {\em spine\/} of $\surf$.
Matveev (see Kudryavtseva~\cite{Kudr}) notes that the space of Morse functions with isotopic spines is path-connected.
He further suggested elementary transformations of spines which induce \sgmhom-homotopies of Morse function and showed that any two spines can be connected by a finite sequence of these transformations.

Sharko~\cite{Sharko} reduced the problem to minimal Morse functions on a surface $\surf$ with only one positive and only one negative boundary component.
Such a surface can be regarded as a ``framed'' chords diagram in which the union of all chords and a negative boundary component constitute the spine of $\surf$.
Notice that $\pi_1\surf$ is free.
Choose a basis of this group.
Then the edges of any other chords diagram in $\surf$ can be written down as words in the terms of a given basis.
These words also form the basis of $\pi_1\surf$ and determine chords diagram up to equivalence.
Moreover, by the well-known Nielsen theorem any two bases of a finitely generated free group are related by a finite sequence of Nielsen transformations.
Sharko proved that Nielsen transformations yield \sgmhom-homotopies between corresponding Morse functions, and that Morse functions with equivalent diagrams are \sgmhom-homotopic.

The extension of the proof of~\cite{Maks:PhD} for $\onemanif=\Circle$ and all surfaces is given in the Appendix.

\subsection{Plan of the present proof.}
First the problem will be reduced to the case when $\gfunc=\func\circ\hdif$, where $\hdif$ is a diffeomorphism of $\surf$ and $\func$ is of a special ``canonical'' form.
It is convenient to say that a diffeomorphism $\hdif$ is $\func$-admissible if $\func\sgmh\func\circ\hdif$. 
Using a special type of $\func$, we will choose system of generators for $\HMbd$ and show that if $\onemanif=\RRR$, then all of them are $\func$-admissible.
This will prove the \MainTheorem\ for this case.

For the case $\onemanif=\Circle$, $\surf$ is orientable, and $\func$ is not null-homotopic we shall see that one of the generators chosen above is not $\func$-admissible.
Nevertheless, since $\func$ and $\func\circ\hdif$ are homotopic,
it will be possible to reduce the problem to the case when $\hdif$ acts trivially on the homology group $H_1(\manif,\partial\manif)$, i.e.\! $\hdif$ belongs to the Torelli group of $\surf$.
Generators of this group are known from~\cite{Powell},~\cite{Johnson},~\cite{Mess}.
This information will allow us to show that $\func\sgmh\func\circ\hdif$.

\subsection{Structure of the paper.}
In Section~\ref{subsect:CuttingAlongLevelSet} we prove some technical results concerning to Morse mappings to the circle.
In Section~\ref{subsect:ReebGraph} we recall the definition of the Kronrod-Reeb graph of a Morse mapping and define ``canonical'' Morse mappings.
In Section~\ref{sect:reduc_probl} we reduce the \MainTheorem\ to the case when $\func$ is canonical and $\gfunc$ differs from $\func$ by a diffeomorphism. This was done by Kudryavtseva in~\cite{Kudr} for Morse functions. We consider the case $\onemanif=\Circle$.
In Section~\ref{sect:admiss} we show that elementary diffeomorphisms generating mapping class group $\HMbd$ of $\surf$
(Dehn twists, boundary and crosscap slides) preserve certain Morse functions.
In Section~\ref{sect:MapClGr_BndSurf} we recall the generators of mapping class groups for surfaces with boundary. 
Every canonical Morse mapping gives a ``canonical'' set of such generators whose admissibility (or nonadmissibility) for this map is almost obvious.
We also complete the \MainTheorem\ for $\onemanif=\RRR$ (statement (i) of Lemma~\ref{lm:canon_gen_for_canon_Morse}).

In Section~\ref{sect:proof_main_th_S1} we give the plan of the proof of \MainTheorem\ for the case $\surf$ is orientable and $\onemanif=\Circle$.
For this in Section~\ref{sect:symplect_grp} we consider the stabilizers of elements of $\ZZZ^{2g}$ with respect to the action of the symplectic groups $\SpZg$, in Section~\ref{sect:minimal_Morse_maps} we study minimal Morse functions.  Section~\ref{sect:minimizat_int_levsets} includes one technical lemma.
Finally in Sections~\ref{sect:prf_th_admis_i}-\ref{sect:prf_lm_reduce_torel} we complete the proof.


\section{Cutting $\surf$ along a regular level-set of $\func$.}
\label{subsect:CuttingAlongLevelSet}
We prove here two lemmas which will be used in the proof of Proposition~\ref{pr:toReebGr}.

Let $\cval$ be a regular value of a Morse mapping $\func:\surf\to\Circle$.
Then $\func^{-1}(\cval)$ is a disjoint union of \SCCs\ on $\surf$.
Suppose that $\func^{-1}(\cval) \cap \partial\surf=\varnothing$.
We cut $\surf$ along $\func^{-1}(\cval)$ and denote the new surface by
$\tsurf = \tsurf\frval{\func}{\cval}$.
Similarly, we cut $\Circle$ at $\func(\cval)$ and obtain $[0,1]$.
Let $\psurf:\tsurf\to\surf$
and $\pcirc:[0,1]\to\Circle$ be the corresponding factor-maps,
where $\pcirc(t) = e^{2\pi i t}$, $t\in[0,1]$.
Then there exists a Morse function $\tfunc:\tsurf\to[0,1]$
such that the following diagram is commutative:
\chEq\begin{equation}\label{equ:CD_corr_sgmhom}
\begin{CD}
\tsurf @>{\tfunc}>> [0,1] \\
@V{\psurf}VV @VV{\pcirc}V \\
\surf @>{\func}>> \Circle.
\end{CD}
\end{equation}
Thus
\chEq\begin{equation}\label{equ:fx_e_tfpx}
\func(x) = \exp \left({2\pi i \tfunc(\psurf^{-1}(x))}\right), \ \forall x\in\surf.
\end{equation}
Denote $\btsurf{0}=\tfunc^{-1}(0)$, $\btsurf{1}=\tfunc^{-1}(1)$,
and $\btsurf{} = \btsurf{0}\cup\btsurf{1}$.
Then there is a natural corresponding between \sgmhom-homotopies $\tfsgm_t$
of $\tfunc$ with respect to some neighborhood of $\btsurf{}$
and \sgmhom-homotopies $\fsgm_t$ of $\func$ with respect to some neighborhood of $\tmpcurve$.
The corresponding maps $\tfsgm_t$ and $\fsgm_t$ are
related by the commutative diagram~\eqref{equ:CD_corr_sgmhom}.

Since $\surf$ is connected, it follows that
every connected component $\comp$ of $\tsurf$ intersects $\btsurf{}$.
However, it is possible that $\comp\cap\btsurf{i}=\varnothing$ for some $i=0,1$.
Thus the components of $\tsurf$ can be divided into the following mutually disjoint sets
\chEq\begin{equation}\label{equ:comp_QQQ}
\compAA=\compAA\frval{\func}{\cval},
\qquad
\compAB=\compAB\frval{\func}{\cval},
\qquad
\compBB=\compBB\frval{\func}{\cval}
\end{equation}
that (respectively) intersect only $\btsurf{0}$,
intersect both sets $\btsurf{0}$ and $\btsurf{1}$,
and intersect only $\btsurf{1}$.

It follows that for every connected component $\comp$ of $\compAB\frval{\func}{\cval}$
and $t\in[0,1]$ we have $\comp\cap\func^{-1}(t)\not=\varnothing$.

\begin{figure}[ht]
\centering
\begin{tabular}{ccc}
\includegraphics[height=2.1cm]{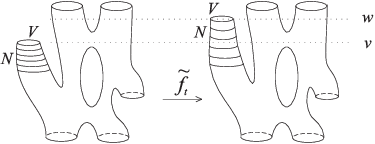} &
\qquad &
\includegraphics[height=2.1cm]{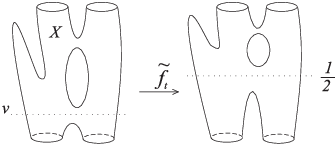} \\
a) & & b)
\end{tabular}
\caption{}
\protect\label{fig:prf12}
\end{figure}

\begin{lem}\label{lm:two_sgm_hom}
{\em(1)}
Let $\bcomp$ be an $\tfunc$-positive (resp. $\tfunc$-negative)
component of $\partial\tsurf$ and $\vval = \tfunc(\bcomp)$.
Then for every $\wval>\vval$ (resp. $\wval<\vval$)
there exists a \sgmhom-homotopy $\tfunc_t$
changing $\tfunc$ only in arbitrary small neighborhood of $\bcomp$
and such that $\tfunc_1(\bcomp)=\wval$, see Figure~\ref{fig:prf12}a).

{\em(2)}
Let $\comp$ be a connected component of $\tsurf$.
For every $\wval\in(0,1)$ there exists a \sgmhom-homotopy
$\tfunc_t:\tsurf\to[0,1]$ such that $\tfunc_0=\tfunc$,
$\tfunc_t = \tfunc$ on $(\tsurf\setminus\comp)\cup\btsurf{}$,
and $\tfunc_1^{-1}(\onehalf)\cap \comp = \tfunc^{-1}(\wval)\cap \comp$,
see Figure~\ref{fig:prf12}b).

{\em(3)}
Let $\comp$ be a connected component of $\tsurf$.
Then there exists a \sgmhom-homotopy $\tfunc_t:\tsurf\to[0,1]$ with $\tfunc_0=\tfunc$ and
$\tfunc_t = \tfunc$ on $(\tsurf\setminus\comp)\cup\btsurf{}$,
such that $\tfunc_1^{-1}(\onehalf)\cap \comp=\varnothing$, whenever $\comp\subset\compAA\cup\compBB$,
and
$\concomp{\tfunc_1^{-1}(\onehalf)\cap \comp}=1$, whenever $\comp\subset\compAB$.
\end{lem}

\chFig
\proof 
(1)
Suppose that $\bcomp$ is an $\tfunc$-positive component of $\partial\tsurf$.
By definition, $\tfunc$ has no critical points on $\bcomp$.
Then there exist an $\eps>0$, a neighborhood $\bNbh$ of $\bcomp$,
and a diffeomorphism $\hdif:\Circle \times (\vval-2\eps, \vval] \to \bNbh$
such that $\hdif(\Circle\times\{\vval\}) = \bcomp$
and $\tfunc\circ\hdif(x,t) = t$ for $(x,t)\in\Circle\times(\vval-2\eps, \vval]$.

Let $\Hisot_t$ be an isotopy of $\RRR$ fixed on $(-\infty,\vval-\eps]$
and such that $\Hisot_1(\vval)=\wval$.
Then the \sgmhom-homotopy $\tfunc_t$ defined by the formulas:
$\tfunc_t(x) = \tfunc(x)$ for $x\in\surf\setminus\bNbh$
and $\tfunc_t(x) = \Hisot_t\circ \tfunc(x)$ for $x\in\bNbh$
satisfies the statement (1) of lemma.
The proof for $\tfunc$-negative components is similar.

(2)
Notice, that for any $\vval\in(0,1)$ there exists an isotopy $\Hisot_t$ of $\Rline$ fixed near $0$ and $1$ and such that $\Hisot_1(s)=\onehalf$.
Then the \sgmhom-homotopy $\tfunc_t:\tsurf\to[0,1]$ defined by the formulas:
$\tfunc_t = \Hisot_t \circ \tfunc$ on $\comp$ and
$\tfunc_t = \tfunc$ on $\tsurf\setminus\comp$
satisfies the statement (2) of lemma.

(3)
It follows from the definition that for every connected component $\comp$ of $\compAA\cup\compBB$ there exists a number $\vval\in(0,1)$ such that $\tfunc^{-1}(\vval)\cap\comp=\varnothing$.
Therefore, if $\comp\subset\compAA\cup\compBB$, then our statement follows from (2).

Let $\comp\subset\compAB$.
If for some $i=0,1$ the intersection $\comp\cap\btsurf{i}$ is connected,
then for every $t$ in some neighborhood of $i$ we have that $\comp\cap\tfunc^{-1}(t)$ is connected.
By (1) of this lemma we can choose $t=\onehalf$.

Suppose now that the intersections $\comp\cap\btsurf{i}$, $i=0,1$ are not connected.
By (1) and (2) we assume that
$$
0 < \tfunc(\psurf^{-1}(\negbndf) \cap \comp) < \oneforth <
\tfunc(\Sigma_{\tfunc} \cap \comp) < \onehalf < \tfunc(\psurf^{-1}(\posbndf) \cap \comp) < 1,
$$
where $\Sigma_{\tfunc}$ is the set of critical points of $\tfunc$.
Thus all critical values of $\tfunc|_{\comp}$ belong to $(\oneforth,\onehalf)$;
the values on $\tfunc$-negative boundary components of $\comp$ except for $\tfunc(\comp\cap\btsurf{0})=0$ are in $(0,\oneforth)$;
and the values on $\tfunc$-positive boundary components of $\comp$ except for $\tfunc(\comp\cap\btsurf{1})=1$ are in $(\onehalf, 1)$.
In particular, $\onehalf$ is a regular value of $\tfunc$.

Denote $n=\concomp{\tfunc^{-1}(\onehalf)}$ and suppose that $n>1$.
Our object is to reduce $n$.
Let $\Fld$ be a gradient-like Morse-Smale vector field of $\comp$ for the function $\tfunc|_{\comp}$.
It follows from Morse theory that the union of $\tfunc|_{\comp}$-positive boundary components $\tposbndfX$ with the set of trajectories that start at saddle critical points of $\tfunc|_{\comp}$ and finish at $\tposbndfX$ is a strong deformation retract of $\comp$.
Since $\comp$ is connected, we see that there exists a saddle critical point $\pnt$ of $\tfunc|_{\comp}$ such that the trajectories starting from
$\pnt$ finish at different components of $\tposbndfX$.
We denote these trajectories by $\omega_1$ and $\omega_2$.

Then (Milnor~\cite{Milnor_h_cob}, Theorem~4.1) there exists a \sgmhom-homotopy $\tfunc_t$ of $\tfunc_0 = \tfunc|_{\comp}$ that changes $\tfunc|_{\comp}$ only in arbitrary small neighborhood of $(\omega_1\cup\omega_2) \cap \tfunc^{-1}(\oneforth, \onehalf]$
such that $\onehalf < \tfunc_1(\pnt) < 1$, but $\tfunc_1(\pnt')<\onehalf$ for all other critical point $\pnt'$ of $\tfunc_1$.
It follows that $\onehalf$ is a regular value for $\tfunc_1$ and the level-set $\tfunc_1^{-1}(\onehalf)$ has precisely $n-1$ connected components.
Now (3) follows by induction on $n$.
\endproof

\begin{lem}\label{lm:throw_out_1_0_comp}
Every Morse mapping $\func:\surf\to\Circle$ is \sgmhom-homotopic to a
Morse mapping $\gfunc$
such that for some regular value $\cval$ of $\gfunc$ we have:
\\ {\em(A)} if $\func$ is null-homotopic, then $\gfunc^{-1}(\cval)=\varnothing$;
\\ {\em(B)} otherwise, $\concomp{\gfunc^{-1}(\cval)}$
is equal to the index of $\func_{\hommap}(H_1(\surf))$ in $H_1(\Circle)$.
\end{lem}
\proof
Let $\cval$ be a regular value of $\func$ such that
$\func^{-1}(\cval)\cap\partial\surf=\varnothing$
and let $n=\concomp{\func^{-1}(\cval)}$.
We cut $\surf$ and obtain the surface $\tsurf=\tsurf\frval{\func}{\cval}$
and the function $\tfunc:\tsurf\to[0,1]$ as above.

By Lemma~\ref{lm:two_sgm_hom}, if
$\compAA\cup\compBB\not=\varnothing$
or if for some connected component $\comp$ of $\compAB$
the intersection $\comp\cap\btsurf{0}$ has more than one component,
then there exists a \sgmhom-homotopy $\tfunc_t$ of $\tfunc$ with respect to
some neighborhood of $\btsurf{}$
such that $\concomp{\tfunc_1^{-1}(\onehalf)}<n$.
As noted above, this \sgmhom-homotopy yields a \sgmhom-homotopy $\func_t$ of $\func=\func_0$ to a Morse mapping $\func_1$ with respect to some neighborhood of $\func^{-1}(\cval)$ such that $\concomp{\func_1^{-1}(\cval_1)}<n$, where $\cval_1 = \pcirc(\onehalf)$ is a regular value of $\func_1$.

Repeating these arguments for $\func_1$ and $\cval_1$, and using induction in $n$ we will obtain a Morse mapping $\func_k$ and its regular value $\cval_k$ such that either (i)~$\func^{-1}_k(\cval_k)=\varnothing$ or
(ii)~$\compAA\frval{\func_k}{\cval_k} = \compBB\frval{\func_k}{\cval_k}=\varnothing$ and for every connected component $\comp$ of $\compAB=\tsurf\frval{\func_k}{\cval_k}$ the intersection $\comp \cap \btsurf{i}\frval{\func_k}{\cval_k}$ is non-empty and connected, whence it is an \SCC.

Suppose that $\func_k$ is null-homotopic.
Then $\func_k$ lifts to a Morse function $\tfunc_k:\surf\to\Rline$ which must have global minimum and maximum.
Therefore, if $\func_k^{-1}(\cval_k)\not=\varnothing$ (case (ii)), then $\compAA\frval{\func_k}{\cval_k}\cup\compBB\frval{\func_k}{\cval_k}\not=\varnothing$, which contradict to (ii).
Hence, $\func^{-1}_k(\cval_k)=\varnothing$.
This proves (A).

Suppose $\func_k$ is not null-homotopic.
For the convenience we denote $\func_k$ by $\func$ and $\cval_k$ by $\cval$.
We will now lift $\func$ onto the covering of $\Circle$ corresponding to the subgroup $\func(H_1(\surf))$ of $H_1(\Circle)$.
Let $m=\concomp{\tsurf}$ and $p_{m}:\Circle\to\Circle$ be the $m$-sheet-covering of $\Circle$ defined by the formula $p_{m}(e^{2\pi i t}) = e^{m 2\pi i t}$, $t\in[0,1]$.

First notice, that the set of connected components of $\tsurf$ admits a natural {\em cyclic\/} ordering.
Indeed, let $\comp_0$ be any component of $\tsurf$.
If $\comp_k$, $(k\geq 0)$ is defined, then there exists a unique connected component $\comp_{k+1}$ of $\tsurf$ such that  $\psurf(\comp_{k+1}\cap\btsurf{0}) = \psurf(\comp_{k}\cap\btsurf{1})$.
Since $\surf$ is connected, it follows that {\em every\/} connected component of $\tsurf$ is numbered in this way.


Then the following formula defines a lifting $\barfunc:\surf\to\Circle$
of $\func$ onto the $m$-sheet covering of $\Circle$: 
$$
\barfunc(x) = \exp{\frac{2\pi i}{m} \left( \tfunc(\psurf^{-1}(x))+k \right)},
\quad  x \in \psurf(\comp_k), \ k=0,\ldots,m-1,
$$
i.e.\! $p_{m}\circ \barfunc = \func$.

Finally, let us prove that the homomorphism $\barfunc_{\hommap}:H_1(\surf) \to H_1(\Circle)$ is onto.
This will imply that index of $\func(H_1(\surf))$ in $H_1(\Circle)$ is $m$.
For every $k=0,\ldots,m-1$ let $\omega_k:[0,1]\to\comp_k$ be a simple arc which is
transversal to level-sets of $\tfunc$ and such that $\tfunc(\omega_k(t))=t$,
$\psurf(\omega_k(1)) = \psurf(\omega_{k+1}(0))$ and
$\psurf(\omega_{m-1}(1)) = \psurf(\omega_{0}(0))$.
Evidently, these arcs constitute an \SCC\ $\omega$ on $\surf$
such that the restriction $\barfunc|_{\omega}$ is homeomorphism, whence $\barfunc_{\hommap}$ is onto.
Thus (B) is proved.
\endproof

\subsection{Orientation of level-sets of $\func$.} \label{sect:or_level_sets}
Suppose that $\surf$ is orientable.
Let $\cval\in\Circle$ be a regular value of a Morse mapping $\func:\surf\to\Circle$, $\lset=\func^{-1}(\cval)$ be the corresponding level-set of $\func$, and $\Fld$ be a gradient vector field for $\func$ taken in some Riemannian metric on $\surf$.
Then the orientation of $\surf$ together with $\Fld$ yields an orientation of $\lset$
so that the homology class of an oriented cycle 
$[\func^{-1}(\cval)]\in H_{1}(\surf,\partial\surf)$ does not depend on a particular choice of a regular value $\cval$ and even on the homotopy class of $\func$. 
For every $x\in\lset$ let $v_{x}$ be a tangent vector to $\lset$ at $x$ such that the pair $(v_{x},\grad\func(x))$ gives a positive orientation of $\surf$. Then the orientation of  $\lset$ defined by $v_x$ satisfies the conditions of the previous sentence.

Let $\xi\in H^{1}(\Circle)$ be a generator that defines the positive orientation of $\Circle$ and $\intform{\cdot}{\cdot}$ be an intersection form on $H_1(\surf,\partial\surf)$.
Then for every oriented \SCC\ $\tmpcurve:\Circle\to\surf$ regarded as an element of $H_1(\surf)$ we have 
\chEq\begin{equation}\label{equ:intersetions}
\func(\xi)(\tmpcurve)=\intrsc{\lset}{\tmpcurve}=\deg(\func|_{\tmpcurve}).
\end{equation}

Since $\func$ is constant on boundary components of $\surf$ and is not null-homotopic it follows that $\func(\xi)\not=0$ in $H^{1}(\surf,\partial\surf)$.
The intersection form $\intform{\cdot}{\cdot}$ on $\surf$ yields an isomorphism 
$\phi:H^{1}(\surf,\partial\surf) \to H_{1}(\surf,\partial\surf)$ which by Eq.~\eqref{equ:intersetions}, maps $\func(\xi)$ to the homology class $[\lset]$.

In particular, if $\hdif:\surf\to\surf$ is a diffeomorphism such that 
$\func\circ\hdif$ and $\func$ are homotopic, then it follows that $\hdif^{\hommap}(\func(\xi))=\func(\xi)$ in $H^{1}(\surf,\partial\surf)$ and $\hdif_{\hommap}([\lset])=[\lset]$ in $H_{1}(\surf,\partial\surf)$.


\section{Kronrod-Reeb graph of a Morse mapping}\label{subsect:ReebGraph}
Let $\func:\surf\to\onemanif$ be a Morse mapping, $\cval\in\onemanif$,
and $\tmpcurve$ be a connected component of $\func^{-1}(\cval)$.
We call $\tmpcurve$ {\em regular\/} if it contains no critical points of $\func$;
otherwise $\tmpcurve$ is {\em critical}.

Consider the partition of $\surf$ by the connected components of level-sets of $\func$.
The factor-space $\Reebf$ of $\surf$ by this partition has the structure of a one-dimensional CW-complex and is called the {\em Kronrod-Reeb graph\/} or \KR-graph of $\func$ (see e.g.~\cite{Kronrod,Kudr,SharkoUMZ}).
There is a unique decomposition
$$ \begin{CD} \func: \surf @>{\pReebf}>> \Reebf @>{\funcR}>> \onemanif, \end{CD} $$
where $\pReebf$ is a factor map and for every {\em open\/} edge $e$ of $\Reebf$ the restriction $\funcR|_{e}$ is a local homeomorphism.
Notice that the orientation of $\onemanif$ yields a unique orientation of $e$ preserved by $\funcR$.
The mapping $\funcR$ will be called {\em \KR-map\/} associated with $\func$.

\begin{figure}[ht]
\centering
\includegraphics[width=7cm]{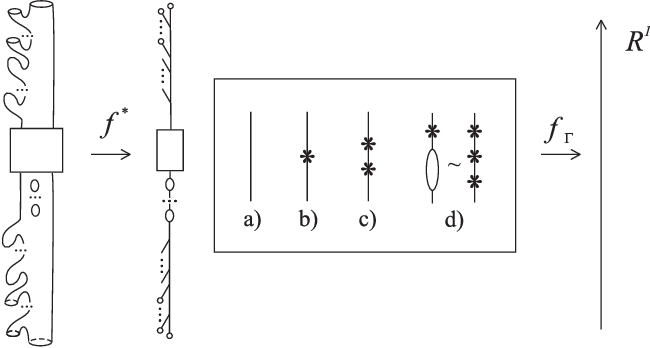} 
\caption{\KR-graphs and \KR-maps of a \canonical\ Morse function $\surf\to\RRR$.}
\protect\label{fig:canon_Reeb_graph_R}
\end{figure}

The vertices of $\Reebf$ correspond to the critical components of level-sets of $\func$ and to the boundary circles of $\surf$.
The last type vertices will be denoted on the \KR-graph by circles $\circ$
(see e.g. Figure~\ref{fig:canon_Reeb_graph_R}).
Notice that for non-orientable surfaces \KR-graphs can possess vertices of degree $2$ (e.g. \cite{Kudr}).
We will denote these vertices by stars $\mystar$.

Let $\func, \gfunc:\surf\to\onemanif$ be Morse mappings.
By {\em isomorphism\/} between their \KR-graphs we will mean a homeomorphism $\Reebg\to\Reebf$ preserving orientations of edges and the sets of $\circ$- and $\mystar$-vertices.

We will say that their \KR-maps $\funcR$ and $\gfuncR$ are {\em equivalent} provided there exist a preserving orientation diffeomorphism $\Cdif$ of $\onemanif$ and an isomorphism $\Rdif:\Reebg\to\Reebf$ such that in the following diagram the right square is commutative:
\chEq\begin{equation}\label{equ:equiv_ReebG}
\begin{CD}
 \surf @>{\pReebg}>> \Reebg @>{\gfuncR}>> \onemanif \\
@V{\hdif}VV @V{\Rdif}VV  @V{\Cdif}VV \\
 \surf @>{\pReebf}>> \Reebf @>{\funcR}>> \onemanif
\end{CD}
\end{equation}

The mappings $\func$ and $\gfunc$ are said to be {\em equivalent} provided there exists a diffeomorphism $\hdif$ of $\surf$ such that $\func\circ\hdif = \Cdif\circ \gfunc$.
In this case there is a unique equivalence $\Rdif$ between \KR-maps of $\func$ and $\gfunc$ such that the whole diagram~\eqref{equ:equiv_ReebG} is commutative.

A Morse mapping $\func$ is called {\em generic\/} if every level-set of $\func$ contains at most one critical point.
Let $\func$ be a generic Morse mapping.
If $\manif$ is orientable, then the degree of each vertex of $\Reebf$ is either $1$ or $3$.
If $\manif$ is non-orientable, then $\Reebf$ can possess vertices of degree $2$.

The following lemma is well-known.
Its different variants can be found in~\cite{BolsFom,Kudr,Kulinich,SharkoUMZ}.
\begin{lem}\label{lm:gen_equiv}
Two generic Morse mappings $\func$ and $\gfunc$ having equivalent \KR-maps are equivalent.
\qed
\end{lem}

\begin{figure}[ht]
\centering
\includegraphics[width=3cm]{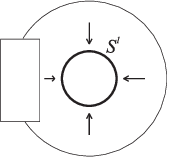} 
\caption{\KR-graphs and \KR-maps of a \canonical\ Morse mapping $\surf\to\Circle$.}
\protect\label{fig:canon_Reeb_graph_S}
\end{figure}

We say that a Morse mapping $\func$ is {\em \canonical\/} if its \KR-map is equivalent to that drawn in Figures~\ref{fig:canon_Reeb_graph_R} or~\ref{fig:canon_Reeb_graph_S}.

First consider the case $\onemanif=\RRR$, see Figure~\ref{fig:canon_Reeb_graph_R}.
The part of \KR-graph under the rectangle corresponds to the following cases of $\surf$:

$a)$ $\surf$ is orientable;

$b)$ $\surf$ is non-orientable of odd genus $\genus$;

$c)$ $\surf$ is non-orientable of even genus $\genus$;

$d)$ $\surf$ is non-orientable, $\genus\geq 3$ and is odd. In this case we will use two types of \canonical\ Morse functions shown in Figure~\ref{fig:canon_Reeb_graph_R}. 
They are related by a \sgmhom-homotopy, see~\cite{Kudr}.

For the case $\onemanif=\Circle$ a \canonical\ Morse mapping $\func:\surf\to\Circle$ can be described as follows: there is a regular value $\cval$ of $\func$ such that $\tmpcurve=\func^{-1}(\cval)$ is a \SCC.
Moreover, if we cut $\surf$ along $\tmpcurve$, then the restriction of $\func:\surf\setminus\tmpcurve\to\Circle\setminus\cval$ is a canonical Morse function.
Its \KR-graph is hidden behind the rectangle, see Figure~\ref{fig:canon_Reeb_graph_S}.

Notice also, that a \canonical\ Morse mapping is generic and the homomorphism $\func_{\hommap}:H_1(\surf)\to H_1(\Circle)$ is onto.

\chFig

\chFig

\begin{lem}\label{lm:canon_mf_equiv}
Let $\func,\gfunc:\manif\to\onemanif$ be two canonical Morse mappings of same critical type $\Crtype{\func}=\Crtype{\gfunc}$.
Then they are equivalent.

Moreover, there is a \sgmhom-homotopy of $\gfunc$ to a canonical Morse mapping $\gfunc_1$ such that $\gfunc_1 = \func \circ \hdif$, where $\hdif$ is a diffeomorphism of $\manif$.
\end{lem}
\proof
Evidently, \KR-graph and \KR-map of a canonical Morse mapping is determined by the numbers $\crpt{0}$, $\crpt{2}$, $b_{+}$, $b_{-}$ and the (orientable or non-orientable) genus $\genus$ of $\surf$.
Notice that $\crpt{1}$ is related to these numbers via Euler characteristic.

Hence the condition $\Crtype{\func}=\Crtype{\gfunc}$ implies that \KR-maps of $\func$ and $\gfunc$ are equivalent.
Then by Lemma~\ref{lm:gen_equiv}, $\func$ and $\gfunc$ are equivalent, i.e. $\pdif\circ\gfunc=\func\circ\hdif$, where $\pdif$ is a preserving orientation diffeomorphism of $\onemanif$ and $\hdif$ is a diffeomorphism of $\manif$.
It follows that $\pdif$ is isotopic to $\id_{\manif}$. Let $\pdif_{t}$ be an isotopy of $\pdif=\pdif_1$ to $\id_{\manif}=\pdif_0$.
Then $\gfunc_{t} = \pdif_{t}\circ\gfunc$ is a \sgmhom-homotopy of $\gfunc=\gfunc_{0}$ to $\gfunc_1=\pdif_{1}\circ\gfunc=\pdif\circ\gfunc=\func\circ\hdif$.
\endproof


\section{Reduction of the problem}\label{sect:reduc_probl}
Let $\func,\gfunc:\manif\to\onemanif$ be two Morse mappings such that
$\Crtype{\func}=\Crtype{\gfunc}$.
We have to prove that $\func\sgmh\gfunc$.

In this section we reduce the proof of \MainTheorem\ to the case
when $\func$ and $\gfunc$ are canonical, and $\gfunc=\func\circ\hdif$, where $\hdif$ is a diffeomorphism of $\surf$.
This was done in~\cite{Kudr} for the case $\onemanif=\RRR$. 
Let $\onemanif=\Circle$.

\subsection{Step 1.}
It may be assumed that the 
homomorphism $\func_{*}=\gfunc_{*}:H_{1}(\manif) \to H_{1}(\Circle)$ is surjective.
In particular, $\func$ and $\gfunc$ are not null-homotopic.
This also implies that $\surf$ is neither a sphere nor a projective plane (with holes if $\partial\surf\not=\varnothing$).

Indeed, suppose that the homomorphism $\func_{\hommap}=\gfunc_{\hommap}$ is not onto.
Let $p:\tCircle\to\Circle$ be the covering of $\Circle$ corresponding to the subgroup
$\func_{\hommap}(H_1(\surf)) \subset H_1(\Circle)=\pi_1(\Circle)$
and $\tfunc,\tgfunc:\surf\to\tCircle$ be some liftings of $\func$ and $\gfunc$ respectively which are evidently Morse.
\begin{lem}\label{lm:reduce_to_H1_onto}
$\func\sgmh\gfunc$ iff $\tfunc\sgmh\tgfunc$. \qed
\end{lem}
The proof is easy and is left to the reader. It can be found in~\cite{Maks:PhD}.

\subsection{Step 2.}
We may assume that $\func$ and $\gfunc$ are canonical due to the following statement:
\begin{prop}[\cite{Kudr}]\label{pr:toReebGr}
Every Morse mapping $\func:\surf\to\onemanif$ such that the homomorphism $\func_{\hommap}(H_1(\surf)) \subset H_1(\Circle)=\pi_1(\Circle)$ is onto is \sgmhom-homotopic to a \canonical\ one. 
\end{prop}
It follows from this proposition that $\func\sgmh\func_1$ and $\gfunc\sgmh\gfunc_1$, where $\func_1$ and $\gfunc_1$ are canonical Morse mappings of same critical type $\Crtype{\func}=\Crtype{\gfunc}$. 
Then by Lemma~\ref{lm:canon_mf_equiv} $\gfunc_1=\func_1\circ\hdif$, where $\hdif$ is a diffeomorphism of $\surf$.
\proof
As noted above, this statement is proved in~\cite{Kudr} (Lemma~10) for closed surfaces and $\onemanif=\RRR$. The proof easily extends to surfaces with boundary.
Suppose that $\onemanif=\Circle$.
Since $\func_{\hommap}$ is onto, it follows from Lemma~\ref{lm:throw_out_1_0_comp}, that $\func$ is \sgmhom-homotopic to a Morse mapping $\func_1$ such that $\acurve=\func_1^{-1}(\cval)$ is an \SCC, where $\cval$ is a regular value of $\func_1$.
Cutting $\surf$ along $\acurve$ as in Section~\ref{subsect:CuttingAlongLevelSet} we obtain a surface $\tsurf$ and a function $\tfunc:\tsurf\to[0,1]$.
Then by the $\RRR$-case of this proposition $\tfunc$ is \sgmhom-homotopic with respect to a neighborhood of $\btsurf{}$ to a canonical Morse function.
This \sgmhom-homotopy yields a \sgmhom-homotopy of $\func$ to a canonical Morse mapping.
\endproof


\section{Admissible diffeomorphisms and curves}\label{sect:admiss}
\begin{defn}
Let $\func:\surf\to\onemanif$ be a Morse mapping.
A diffeomorphism $\hdif:\surf\to\surf$ will be called {\em $\func$-admissible\/} provided  $\func\circ\hdif$ is \sgmhom-homotopic to $\func$. 
Notice that $\func$-admissibility implies that $\hdif$ preserves the sets of $\func$-positive and $\func$-negative components of $\partial\manif$ and that $\func$ and $\func\circ\hdif$ are homotopic.
\end{defn}

Let $\fAdmDif \subset\Diff\surf$ be the set of all $\func$-admissible diffeomorphisms, $\DiffId\surf$ be the identity com\-po\-nent of $\Diff\surf$, and $\PathComp{\func}$ be the path-com\-po\-nent of $\func$ in $\MrsMP$.

\begin{lem}\label{lm:admprop}
$\fAdmDif$ is a group consisting of full isotopy classes, i.e. $\DiffId\surf \subset\fAdmDif$.
Moreover, if \ $\gfunc\in\PathComp{\func}$, then $\gAdmDif = \fAdmDif$.
\end{lem}
\proof
Suppose that $\adif, \bdif \in \fAdmDif$ and let $\func \sgmdef{\Asgm_t} \func\circ\adif$ and $\func \sgmdef{\Bsgm_t} \func\circ\bdif$ be  \sgmhom-homotopies.
Then $\adif\circ\bdif$ and $\adif^{-1}\in\fAdmDif$.
Indeed,
$$ \func  \sgmdef{\Bsgm_t} \func\circ\bdif  \sgmdef{\Asgm_t \circ \bdif}  \func \circ\adif \circ \bdif
\quad\text{and}\quad
\func = \func\circ\adif\circ \adif^{-1} \sgmdef{\Asgm_{1-t}\circ \adif^{-1}} \func\circ\adif^{-1}.$$
Thus $\fAdmDif$ is a group.

If $\adif \sgmdef{\Hisot_t} \adif_1$ is an isotopy, then the homotopy $\func \sgmdef{\Asgm_t \circ\Hisot_t} \func \circ \adif_1$ is a \sgmhom-homotopy. 
Thus $\fAdmDif$ consists of full isotopy classes.

Finally, if $\func \sgmdef{\Bsgm_t} \gfunc$ is a \sgmhom-homotopy,
then $\gfunc \sgmdef{\Bsgm_t} \func
 \sgmdef{\Asgm_t} \func\circ\adif \sgmdef{\Bsgm_{1-t} \circ \adif}  \gfunc\circ \adif.$
Hence $\adif\in\gAdmDif$, i.e. $\fAdmDif \subset \gAdmDif$.
Similarly $\gAdmDif \subset \fAdmDif$.
\endproof

We will now consider three types of ``elementary diffeomorphisms'' and show that they preserve certain simple Morse functions.
\subsection{Dehn twists}
\begin{figure}[ht]
\centering
\begin{tabular}{ccccc}
\includegraphics[height=2cm]{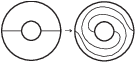} & \qquad \quad &
\includegraphics[height=2cm]{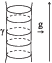} \\
a) & & b) 
\end{tabular}
\caption{Dehn twist}
\label{fig:dehn}
\end{figure}

Let $\tmpcurve$ be a two-sided oriented  \SCC\ in $\surf$.
For definition of a Dehn twist along $\tmpcurve$ see e.g.~\cite{Dehn,Lickor_orient}.
This diffeomorphism is supported in some neighborhood of $\tmpcurve$ and its effect on such a neighborhood is shown in Figure~\ref{fig:dehn}a).

\chFig

\begin{defn}\label{defn:adm_curv_diff}
Let $\tmpcurve$ be a two-sided \SCC\ in $\surf\setminus\partial\surf$.
We say that $\tmpcurve$ is {\em $\func$-admissible} if $\func$ is \sgmhom-homotopic to a Morse mapping $\gfunc$ such that $\tmpcurve$ is a connected component of a regular level-set of $\gfunc$.
\end{defn}

\begin{lem}\label{lm:f_ft}
Let $\tmpcurve \subset \Int\surf$ be a $\func$-admissible oriented \SCC\ in $\surf$.
Then a Dehn twist $\Dtwcurv$ along $\tmpcurve$ is $\func$-admissible.
\end{lem}
\proof
Let $\func \sgmdef{F} \gfunc$ be a \sgmhom-homotopy such that
$\tmpcurve$ is a connected component of a regular level-set of $\gfunc$.
We will construct a Dehn twist $\Dtwcurv$ along $\tmpcurve$ such that $\gfunc=\gfunc\circ\Dtwcurv$.
Then $\Dtwcurv$ is $\gfunc$-admissible, whence by (1) of Lemma~\ref{lm:admprop} $\Dtwcurv$ is also $\func$-admissible.

Since $\tmpcurve$ is a regular component of a level set of $\gfunc$, then there is a regular neighborhood of $\tmpcurve$ which is diffeomorphic to $\Circle\times I$ and such that the function $\gfunc$ is the projection to $I$, see Figure~\ref{fig:dehn}b).
Then there is a Dehn twist $\Dtwcurv$ along $\tmpcurve$ that preserves the sets of the form $\Circle\times\{t\}$. They are level-sets of $\gfunc$, whence $\Dtwcurv$ preserves $\gfunc$.
\endproof

\subsection{Boundary slides}
Let $A$ be an annulus and $C_0,C_1$ be the connected components of $\partial A$.
Divide $C_1$ into four arcs of equal length $l_1,\ldots,l_4$ so that $l_1$ is opposite to $l_3$ and $l_2$ to $l_4$.
Let us identify the opposite points of $l_1$ and $l_3$.
Then the quotient is a \Mobius\ strip $\Mband$ with the hole $C'_1=l_2\cup l_4$. 

\begin{figure}[ht]
\centering
\begin{tabular}{ccccc}
\includegraphics[height=2.5cm]{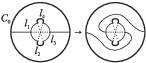} & \qquad \quad &
\includegraphics[height=2.5cm]{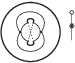} \\
a) & & b) 
\end{tabular}
\caption{Boundary slide}
\label{fig:bslide}
\end{figure}

Let $\tau: A \to A$ be a half-Dehn twist along $C_1$, which exchanges $l_1$ with $l_3$ and $l_2$ with $l_4$ and is identity near $C_0$.
Then $\tau$ yields a certain diffeomorphism $\bsldif$ of $\Mband$ that ``rotates $C'_1$ by $\pi$ and fixes $C_0$'', see Figure~\ref{fig:bslide}a).

Suppose that $\Mband$ is embedded to $\surf$ so that $C_1$ is mapped onto a connected component $C$ of $\partial\surf$.
Then $\bsldif$ extends by the identity on all of $\surf$.
This diffeomorphism is called a {\em boundary slide} of $C$ along $\Mband$.
 
Notice that our description of boundary slide differs from ones given in~\cite{Korkmaz, Szepietowski}. 
The advantage is an evidence of the symmetry of $\bsldif$.

Now it is easy to see that there is a Morse function $\func:\Mband'\to[0,1]$ having a unique critical point of index $1$ and such that $\func^{-1}(0)=C_0$, $\func^{-1}(1)=C'_1$.
Its critical level sets and the \KR-graph are shown in Figure~\ref{fig:bslide}b).

\chFig

The following lemma is obvious.
\begin{lem}\label{lm:pres_bslide}
$\func:\manif\to\onemanif$ be a Morse mapping on a non-orientable surface $\surf$.
Suppose that \KR-graph of $\func$ has an edge $\edgeR$ such that one of its vertices $v_1$ has degree $2$ and another one $v_2$ corresponds to the boundary component of $\surf$, see Figure~\ref{fig:bslide}b).
Let $N$ be a neighborhood of $\edgeR$ containing no vertices of $\Reebf$ but $\partial\edgeR$.
Then $\Mband=\funcR^{-1}(N)\subset\manif$ is a \Mobius\ band with hole 
and there exists a boundary slide $\bsldif:\surf\to\surf$ of $\funcR^{-1}(v_2)$ along $\Mband$ such that $\func\circ\ydif=\func$.
\qed
\end{lem}

\subsection{Crosscap slides.}\label{sect:crosscap_slide_defn}
This type of diffeomorphisms was introduced by W.~B.~R.~Lickorish~\cite{Lickorish_non_or} and called an $Y$-diffeomorphism. In~\cite{Korkmaz, Szepietowski} the term {\em crosscap slide} is used.
We recall the definition of this diffeomorphism (given in~\cite{BirmanChill}) via oriented double coverings.

\begin{figure}[ht]
\centering
\begin{tabular}{ccccc}
\includegraphics[height=3cm]{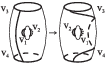} & \qquad \quad &
\includegraphics[height=3cm]{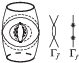} \\
a) & & b) 
\end{tabular}
\caption{Crosscap slide on the orientable covering}
\label{fig:ydif}
\end{figure}

Let $\Kleinbtl$ be a Klein bottle with two holes and
$p:\Torhol\to\Kleinbtl$ be its oriented double covering, where
$\Torhol$ is a torus with $4$ holes.
We can assume that $\Torhol$ is embedded in $\RRR^{3}$ so that it
is symmetrical with respect to the origin $0$. In other words it is invariant under the involution $\xi(x,y,z)=(-x,-y,-z)$ of $\RRR^{3}$, see Figure~\ref{fig:ydif}a).

Let $V_1,\ldots,V_4$ be the connected components of
$\partial\Torhol$ numbered so that $\xi(V_1)=V_2$ and
$\xi(V_3)=V_4$.
Then there is a diffeomorphism
$\tydif:\Torhol\to\Torhol$ which is fixed near $V_3\cup V_4$,
coincides with $\xi$ near $V_1\cup V_2$ and such that $\tydif\circ\xi=\xi\circ\tydif$.
Thus $\ydif$ can be described as a ''rotation'' of $\Torhol$
with respect to $z$-axis by $\pi$ with fixed boundary components $V_3$ and $V_4$.
For example, in Figure~\ref{fig:ydif}a)
an arc and its image under $\xi$ are shown.
It follows that $\tydif$ induces some diffeomorphism $\ydif$ of
$\Kleinbtl$ fixed near $\partial\Kleinbtl$. 

Suppose that $\Kleinbtl\subset\surf$ is embedded in $\surf$.
Then $\ydif$ extends by the identity to a diffeomorphism of $\surf$.
Such a diffeomorphism of $\surf$ is called {\em $Y$-diffeomorphism\/} or {\em crosscap slide\/} based in $\Kleinbtl$.

Notice that there is a Morse function $\tfunc:\Torhol\to\RRR$ with $4$ critical points such that $\tfunc\circ\tydif=\tfunc$, see Figure~\ref{fig:ydif}a), where the critical level-sets of $\tfunc$ are shown.
Then $\tfunc$ yields a unique Morse function $\func:\Kleinbtl\to\RRR$ having $2$ critical points and such that $\func\circ\ydif=\func$.
The \KR-graphs $\Gamma_{\tfunc}$ and $\Gamma_{\func}$ of $\tfunc$ and $\func$ are shown in Figure~\ref{fig:ydif}b).

\chFig

\begin{lem}\label{lm:pres_ydif}
Let $\func:\manif\to\onemanif$ be a Morse mapping on a non-orientable surface $\surf$.
Suppose that \KR-graph of $\func$ has an edge $\edgeR$ with vertices of degree $2$.
Let $N$ be a neighborhood of $\edgeR$ containing no vertices of $\Reebf$ but $\partial\edgeR$.
Then $\Kleinbtl=\funcR^{-1}(N)\subset\manif$ is a Klein bottle with two holes
and there exists an $Y$-diffeomorphism $\ydif:\surf\to\surf$ based in $\Kleinbtl$ such that $\func\circ\ydif=\func$.
\qed
\end{lem}


\section{Mapping class group of a surface with boundary}\label{sect:MapClGr_BndSurf}
Let $\ubdsurf$ be a closed connected surface and $X = \{x_{1},\ldots,x_{n}\}$ be a set of mutually distinct points of $\ubdsurf$.
The {\em extended mapping class group\/} $\HomGr{\surf}{n}$ of $\surf$ is defined to be the group of isotopy classes of diffeomorphisms of $\ubdsurf$ which take $X$ to itself.
The {\em pure extended mapping class group\/} $\PHomGr{\surf}{n}$ of $\surf$ is the group of isotopy classes of diffeomorphisms of $\ubdsurf$ which take $X$ point-wise.
The groups $\HomGr{\ubdsurf}{0}$ and $\PHomGr{\ubdsurf}{0}$ will be denoted by $\HomGr{\ubdsurf}{}$ and $\PHomGr{\ubdsurf}{}$ respectively.

Let $\surf$ be a connected surface with boundary $\partial\surf$ consisting of $n$ connected components $V_1,\ldots,V_n$.
Regarding these components as punctures, we can identify the groups $\HomGr{\surf}{}$ and $\PHomGr{\surf}{}$ with $\HomGr{\ubdsurf}{n}$ and $\PHomGr{\ubdsurf}{n}$.

We recall the sets of generators of $\HomGr{\surf}{}$ and $\PHomGr{\surf}{}$ given in~\cite{Birman,Gervais} for orientable surfaces and in~\cite{Korkmaz} for nonorientable ones.

\subsection{Orientable case.}
Suppose that $\surf$ is orientable.
Consider the following $3$ types of diffeomorphisms of $\bdsurf$:

(1) Let $\DifOr$ be a reversing orientation diffeomorphism of $\bdsurf$.

\begin{figure}[ht]
\centering
\begin{tabular}{ccccc}
\includegraphics[height=2.5cm]{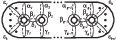} & \qquad \quad &
\includegraphics[height=2.5cm]{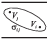} \\
a) & & b) 
\end{tabular}
\caption{The configuration $\CC$. Orientable case.}
\label{fig:generators_or}
\end{figure}

(2) Let $\acurve_i, \bcurve_i, \ccurve_i, \dcurve_i, \ecurve_i$ be the \SCC\ shown in Figures~\ref{fig:generators_or}a), where the bold points denote connected components of $\partial\surf$ divided into two parts (positive and negative components).
We will refer them as \SCCs\ of configuration $\CC$.
Denote by $\Dtw{\acurve_i}$, $\Dtw{\bcurve_i}$, $\Dtw{\ccurve_i}$, $\Dtw{\dcurve_j}$,
$\Dtw{\ecurve_i}$ the corresponding Dehn twists.

(3) For every pair $i<j = 1,\ldots,n$ let $\scurve_{ij}$ be a \SCC\ that separates $\bdsurf$ into two connected components so that one of which is a sphere $S$ with $3$ holes whose boundary components are $\scurve_{ij}$ and the connected components $V_{i}$ and $V_{j}$ of $\partial\surf$, see Figure~\ref{fig:generators_or}b).
Let $\bndif_{ij}$ be a diffeomorphism of $\bdsurf$ with support in $S$ which permutes
boundary components $V_i$ and $V_j$ and preserves all others. 
Evidently, $\bndif_{ij}^2$ is a Dehn twist $\Dtw{\scurve_{ij}}$ along $\scurve_{ij}$.

\begin{thm}[\cite{Birman,Gervais}]\label{thm:gen_or}
The group $\HMbd$ is generated by 
\begin{itemize}
\item[(i)]
$\{\DifOr, \bndif_{ij}: i,j = 1,\ldots,n\}$ if $\genus = 0$; 
\item[(ii)]
$\{\Dtw{l}, \DifOr, \bndif_{ij}: l\in\CC, i,j = 1,\ldots,n\}$ if $g\geq 1$.
\end{itemize}

The group $\PHMbd$ is generated by 
\begin{itemize}
\item[(i)]
$\{\DifOr, \bndif^2_{ij}=\Dtw{\scurve_{ij}}: i,j = 1,\ldots,n\}$ if $\genus = 0$; 
\item[(ii)]
$\{\Dtw{l}, \DifOr: l\in\CC, i,j = 1,\ldots,n\}$ if $g\geq 1$.
\end{itemize}
\end{thm}

\chFig

\subsection{Generators for $\HMbd$. Non-orientable case.}\label{subsect:gen_nor}
\begin{figure}[ht]
\centering
\begin{tabular}{c}
\includegraphics[height=3cm]{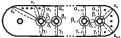} \\
\includegraphics[height=3cm]{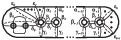}
\end{tabular}
\caption{The configuration $\CC$ for $\genus=2r+1$ and $\genus=2r+2$. Non-orientable case}
\label{fig:generators_nor}
\end{figure}

Suppose that $\surf$ is non-orientable of genus $\genus$, see Figure~\ref{fig:generators_nor}, where the interiors of the shaded disks are removed and then the antipodal points on each boundary component are to be identified.

Consider the following $4$ types of diffeomorphisms of $\bdsurf$:

(1) Let $\ydif$ be a crosscap slide of $\surf$. If $\genus\geq 3$, then we additionally assume that $\ydif^2$ is a Dehn twist along a two-sided separating \SCC\ both components of whose complement are non-orientable.

(2) and (3) 
Similarly to the oriented case we define the configuration $\CC$ of \SCCs\ 
$\acurve_i, \bcurve_i, \ccurve_i, \dcurve_i, \ecurve_i$
shown in Figure~\ref{fig:generators_nor},
\SCCs\ $\scurve_{ij}$, the corresponding Dehn twists and diffeomorphisms $\bndif_{ij}$.

\begin{figure}[ht]
\centering
\begin{tabular}{ccc}
\includegraphics[height=2cm]{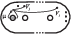} & \quad &
\includegraphics[height=2cm]{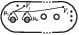}
\end{tabular}
\caption{Boundary slides for $\genus=2r+1$ and $\genus=2r+2$.}
\label{fig:generators_nor_bslides}
\end{figure}

(4) Let $\bsldif_i$ denotes the boundary slide obtaining by sliding the boundary 
component $V_i$ along the loop $\mu$ if $\genus$ is odd and along $\mu_1$ if $\genus$ is even, see Figure~\ref{fig:generators_nor_bslides}.
Also if $\genus$ is even, denote by $\obsldif_i$ the boundary slide obtaining by sliding $V_i$ once along the loop $\mu_2$.

\begin{thm}[\cite{Korkmaz}]\label{thm:gen_nor}
The group $\HMbd$ is generated by 
\begin{itemize}
\item[(i)]
 $\{\bsldif_k, \bndif_{ij}: i,j,k = 1,\ldots,n, i<j\}$ if $\genus=1$;
\item[(ii)]
 $\{\Dtw{\bcurve_{0}}, \ydif, \bsldif_k, \bndif_{ij}: i,j,k = 1,\ldots,n, i<j\}$ if $\genus=2$;
\item[(iii)]
 $\{\Dtw{l}, \ydif, \bsldif_k, \bndif_{ij}: l\in\CC,i,j,k = 1,\ldots,n, i<j\}$  
 if $\genus \geq 3$ is odd;
\item[(iv)]
 $\{\Dtw{l}, \ydif, \bsldif_k, \obsldif_k, \bndif_{ij}: l\in\CC,i,j,k = 1,\ldots,n, i<j\}$  
 if $\genus \geq 4$ is even.
\end{itemize}
Replacing every $\bndif_{ij}$ by $\bndif^2_{ij}=\Dtw{\scurve_{ij}}$ we obtain generators for $\PHMbd$.
\end{thm}

\chFig

\chFig

\subsection{Generators of $\HMbd$ for \canonical\ Morse mapping.}
Given a Morse mapping $\func$, denote by $\HMf$ the subgroup of $\HM$ consisting of diffeomorphisms that preserve the sets of $\func$-positive and $\func$-negative components of $\partial\surf$. Evidently, $\fAdmDif\subset\HMf$.

\begin{lem}\label{lm:canon_gen_for_canon_Morse}
Let $\func:\surf\to\onemanif$ be a canonical Morse mapping. In the case $\onemanif=\Circle$ assume that $\surf$ is orientable.
Then there is a ``canonical'' set of generators for $\HMf$ such that 
\begin{itemize}
\item[(i)]
for the case $\onemanif=\Rline$  all of them are $\func$-admissible, i.e.\!
$\fAdmDif=\HMf$, whence \MainTheorem\ holds for this case;
\item[(ii)]
for $\onemanif=\Circle$ (and orientable $\surf$) all but one of them are also $\func$-admissible.
\end{itemize}
\end{lem}

\begin{rem}\label{rem:partcase}
Recall that we do not give the proof of \MainTheorem\ (by the new method) for the case $\surf$ is non-orientable and $\onemanif=\Circle$.
Therefore we also do not consider this case in Lemma~\ref{lm:canon_gen_for_canon_Morse} since it is more complicated and due to the length of the paper, see also the last paragraph of this section.
\end{rem}

\proof
Let $\func$ be a canonical Morse mapping.
We will construct a set of generators for $\HMbd$ described in Theorems~\ref{thm:gen_or} and~\ref{thm:gen_nor} such that their $\func$-admissibility is rather evident.

First suppose $\surf$ is orientable and embedded in $\RRR^{3}$ as it is shown in Figure~\ref{fig:canon_Reeb_graph_R}. 
Then the canonical Morse mapping $\func$ is just the projection onto the vertical line.

(1) Let $\DifOr$ be a diffeomorphism of $\surf$ that is a symmetry with respect to the plane of this sheet. Then $\DifOr$ reverses orientation of $\surf$ and preserves $\func$, i.e.\! $\func=\func\circ\DifOr$. Thus $\DifOr$ is $\func$-admissible.

\begin{figure}[ht]
\centering
\begin{tabular}{ccccccc}
\includegraphics[height=2.1cm]{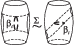} & \quad &
\includegraphics[height=2.1cm]{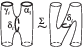}  & \quad & 
\includegraphics[height=2.1cm]{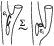}  
\end{tabular}
\caption{$\func$-admissibility of configuration $\CC$}
\label{fig:admis_bde}
\end{figure}

(2) 
Comparing Figures~\ref{fig:canon_Reeb_graph_R} and~\ref{fig:generators_or} we see that $\acurve_i$ and $\ccurve_i$ are regular components of regular level-sets of $\func$, whence the Dehn twists $\Dtw{\acurve_i}$ and $\Dtw{\ccurve_i}$ are admissible.
In Figure~\ref{fig:admis_bde} an $\func$-admissibility of twists $\Dtw{\bcurve_i}$, $\Dtw{\dcurve_i}$ and $\Dtw{\ecurve_i}$ is shown.

\chFig

\begin{figure}[ht]
\centering
\begin{tabular}{ccccccccc}
\includegraphics[height=2.1cm]{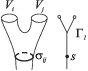} & \qquad &
\includegraphics[height=2.1cm]{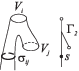} & \qquad & 
\includegraphics[height=2.1cm]{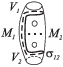}  \\
a) & & b) & & c)
\end{tabular}
\caption{$\func$-admissibility of $\bndif_{ij}$ and $\scurve_{ij}$}
\label{fig:perm_comps}
\end{figure}

(3) 
Let $V_i$ and $V_j$ be two $\func$-positive components of $\partial\surf$.
Then $\func$ is \sgmhom-homotopic to a Morse mapping $\func_{1}$ such that \KR-graph $\ReebGr{\func_{1}}$ of $\func_1$ includes a subgraph $\KrReeb_1$ shown in Figure~\ref{fig:perm_comps}a).
Let $\scurve_{ij}$ be a \SCC\ corresponding to a point $s\in\KrReeb_1$.
Then there exists a diffeomorphism $\bndif_{ij}$ of $\surf_1$ that exchanges $V_i$ and $V_j$, preserves $\func_1$ and $\bndif_{ij}^2$ is a Dehn twist along $\scurve_{ij}$.
Then $\bndif_{ij}$ and $\scurve_{ij}$ are $\func$-admissible.

Now let $V_i$ be $\func$-positive and $V_j$ be $\func$-negative.
In this case a diffeomorphism $\bndif_{ij}$ permuting $V_i$ and $V_j$ is not $\func$-admissible, since it does not preserve the sets of $\func$-positive and $\func$-negative boundary components.
Nevertheless we will now show that its square $\bndif_{ij}^2=\Dtw{\scurve_{ij}}$ is $\func$-admissible.
Consider two cases.

(a) Suppose that $\func$ has at least one critical point of index either $0$ or $2$ or a boundary component different from $V_i$ and $V_j$.
Then $\func$ is \sgmhom-homotopic to a Morse mapping $\func_1$ whose \KR-graph $\ReebGr{\func_{2}}$ includes a subgraph $\KrReeb_2$ shown in Figure~\ref{fig:perm_comps}b).
Then we define $\scurve_{ij}$ to be a \SCC\ corresponding to a point $s\in\KrReeb_2$.
Hence $\scurve_{ij}$ is $\func$-admissible.

(b) Otherwise, $\func$ has no local extremes and $\partial\surf=V_1\cup V_2$. Let $\scurve_{12}$ be a \SCC\ that intersects every $\ccurve_{i}$ but no other \SCCs\ of configuration $\CC$, separates $\surf$ in two components $\surf_1$ and $\surf_2$ such that $\surf_1$ is disk with two holes $V_1$ and $V_2$, see Figure~\ref{fig:perm_comps}c).

We claim that $\scurve_{12}$ is not $\func$-admissible.
Otherwise the restriction of $\func$ to $\surf_{2}$ must have extremes, which could be taken only on boundary components different from $V_1$ and $V_2$ or at critical points of indexes $0$ and $2$. But all of them are absent on $\surf_{2}$.

Nevertheless, it is well-known that a Dehn twist $\Dtw{\scurve_{12}}$ is a product of Dehn twists along \SCCs\ of configuration $\CC$ except for $\ccurve_{i}$.
Hence a Dehn twist $\Dtw{\scurve_{12}}$ is $\func$-admissible.

\chFig

Suppose that $\surf$ is non-orientable of genus $\genus$ (see Figure~\ref{fig:generators_nor}) and let $\func$ be a \canonical\ Morse mapping as in Figure~\ref{fig:canon_Reeb_graph_R}.
Again we define the generators of $\HMbd$ associated with $\func$.

(1) For the case $\genus\geq 2$ we will now define an $\func$-admissible crosscap slide.
If $\genus$ is odd then, $\Reebf$ has an edge $\edgeR$ with vertices of degree $2$.
Otherwise, $\func$ is \sgmhom-homotopic to a Morse function $\func_1$ whose \KR-graph has such an edge, see Figure~\ref{fig:canon_Reeb_graph_R}d).
Then by Lemma~\ref{lm:pres_ydif}, there exists a crosscap slide $\ydif$ such that $\func=\func\circ\ydif$ or $\func_1=\func_1\circ\ydif$ in the second case.
Hence $\ydif$ is $\func$-admissible.

Definition and $\func$-admissibility of generators of types (2) and (3) are similar to the orientable case.
We need to verify the admissibility of $\bcurve_{0}$ and $\dcurve_{0}$ for the case $\genus=2r\geq 2$.

\begin{figure}[ht]
\centering
\begin{tabular}{ccccc}
\includegraphics[height=2.5cm]{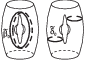} & \quad &
\includegraphics[height=2.5cm]{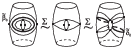}   \\
a) & & b) 
\end{tabular}
\caption{$\func$-admissibility of $\bcurve_{0}$ and $\dcurve_{0}$}
\label{fig:admis_b0d0}
\end{figure}

Let $N$ be a neighborhood of $\edgeR$ defined just above containing no vertices of $\Reebf$ but $\partial\edgeR$.
Then $\Kleinbtl=\funcR^{-1}(N)\subset\manif$ is a Klein bottle with two holes.
Let $p:\Torhol\to\Kleinbtl$. Then $\Torhol$ is a torus with four holes.
We can assume that the function $\tfunc=\func\circ p:\Torhol\to\RRR$ coincides with 
one defined in Section~\ref{sect:crosscap_slide_defn}, see Figure~\ref{fig:ydif}.
Since $\bcurve_{0}$ and $\dcurve_{0}$ are two sided, their inverse images $\tbcurve_{0}=p^{-1}(\bcurve_{0})$ and $\tdcurve_{0}=p^{-1}(\dcurve_{0})$ in $\Torhol$ consist of pair of disjoint \SCC. They are shown in Figure~\ref{fig:admis_b0d0}a). 
 
It is shown in Figure~\ref{fig:admis_b0d0}b) that $\tbcurve_{0}$ is a regular level-set of $\tfunc$. 
This figure also shows a symmetrical \sgmhom-homotopy of $\tfunc$ fixed near $\partial\Torhol$ which makes $\tdcurve_{0}$ a regular level-set. 
Hence $\tbcurve_{0}$ and $\tdcurve_{0}$ are $\tfunc$-admissible, whence 
$\bcurve_{0}$ and $\dcurve_{0}$ are $\func$-admissible.

\chFig

(4) It remains to construct $\func$-admissible boundary slides $\bsldif_{i}$ and $\obsldif_{i}$.
Let $V_i$ be a connected component of $\partial\surf$ and $z_i\in\Reebf$ be the corresponding $\circ$-vertex.

\begin{figure}[ht]
\centering
\begin{tabular}{c}
\includegraphics[height=1cm]{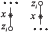}\\
\end{tabular}
\caption{} 
\label{fig:admis_blsides}
\end{figure}

First suppose that $\genus$ is odd, so $\Reebf$ has a unique vertex $x$ of degree $2$.
Then $\func$ is \sgmhom-homotopic to a Morse function $\func_1$ such that $z_i$ and $x$ will be the vertices of same edge, see Figure~\ref{fig:admis_blsides} for the cases when $z_i$ is $\func$-negative or $\func$-positive.
Then by Lemma~\ref{lm:pres_bslide}, there exists a boundary slide $\bsldif_{i}$ of $V_i$ preserves $\func_1$. Whence $\bsldif_{i}$ is $\func$-admissible.

If $\genus$ is even, then $\Reebf$ has two vertices $x_1$ and $x_2$ of degrees $2$.
As in the previous case we define $\func$-admissible boundary slices $\bsldif_{i}$ for $V_i$ and $x_1$, and $\obsldif_{i}$ for $V_i$ and $x_2$.

\chFig

Consider now the case $\onemanif=\Circle$.
Let $\cval\in\Circle$ be a regular value of $\func$ and $\acurve_1=\func^{-1}(\cval)$ such that the restriction of $\func$ to $\surf\setminus\acurve_1$ is a canonical Morse function to $\Circle\setminus\cval$.

\begin{figure}[ht]
\centering
\begin{tabular}{c}
\includegraphics[height=4cm]{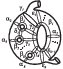}\\
\end{tabular}
\caption{Configuration $\CC$ if $\surf$ is orientable and $\onemanif=\Circle$}
\label{fig:canon_S1_gen}
\end{figure}

Suppose that $\surf$ is orientable. 
Then the definition configuration of the $\CC$ associated with $\func$ is shown in Figure~\ref{fig:canon_S1_gen}, where $\func$ is the ``projection'' to $\bcurve_{1}$.
Similarly to the previous case we can define a diffeomorphism $\DifOr$, Dehn twists along the  \SCCs\ of configuration $\CC$, and permutations of boundary components $\bndif_{ij}$.
The same arguments as in the case $\onemanif=\RRR$ show that all of them are admissible, except for $\bcurve_1$, since $\func$ and $\func\circ\Dtw{\bcurve_1}$ are not even homotopic.

If $\surf$ is non-orientable, then the surface $\surf\setminus\acurve_1$ can be orientable or non-orientable as well.
We do not consider this case, see~\ref{rem:partcase}.
\endproof

\chFig


\section{Proof of \MainTheorem.}\label{sect:proof_main_th_S1}
The case $\onemanif=\RRR$ is proved in statement (i) of Lemma~\ref{lm:canon_gen_for_canon_Morse}.
Before processing with the case $\onemanif=\Circle$ we recall the definition of Torelli group and its generators.

\subsection{Torelli group $\TrlGrM$.}
Let $\surf$ be a closed orientable surface.
Then the Torelli group of $\surf$ is a subgroup $\TrlGrM$ of $\PHMbd=\HMbd$ consisting of diffeomorphism of $\surf$ trivially acting on the homology group $H_1(\surf)$.
Evidently, $\TrlGrM$ is a normal subgroup in $\PHMbd$.

Suppose now that $\partial\surf\not=\varnothing$.
Let us glue every connected component of $\partial\surf$ by $2$-disk and denote the obtained closed surface by $\ubdsurf$.
Then we obtain an epimorphism $j: \PHomGr{\surf}{} \to \PHomGr{\ubdsurf}{}$ induced by the inclusion $\surf\subset\ubdsurf$, see~\cite{Birman}.
Define the Torelli group $\TrlGrM \subset \PHomGr{\surf}{}$ of $\surf$ to be the inverse image $j^{-1}(\TT(\ubdsurf))$.

The following theorem describes the generators of $\ker j$.
\begin{thm}[\cite{Birman:BraidGroup, Birman}]\label{th:gen_of_ker_j}
Let $\acurve_i$ and $\bcurve_i$ be the curves of configuration $\CC$ on $\surf$.
For every component $V_j$ of $\partial\surf$ let $\acurve_{ik}$ ($\bcurve_{ik}$) be an \SCC\ which together with $\acurve_i$ ($\bcurve_{i}$) bounds in $\surf$ a cylinder with a hole $V_i$.
Then the kernel of $j$ is generated by the following diffeomorphisms: $s_{ik}=\acurve_i\circ\acurve_{ik}^{1}$ and
$r_{ik}=\bcurve_i\circ\bcurve_{ik}^{1}$.
\end{thm}

\begin{thm}[\cite{Birman:SiegelModGr, Powell, Johnson, Mess}]\label{th:gen_Torel_grp}
The Torelli group $\TrlGrM$ of $\surf$ is generated by the following types of diffeomorphisms:

{\rm (a)} Dehn twists along \SCC\ separating $\surf$ (if $\genus=2$ then these diffeomorphisms generate all the group $\TrlGrM$,~\cite{Mess});

{\rm (b)} products of Dehn twists of the form $\Dtw{\tmpcurve_1}\circ\Dtw{\tmpcurve_2}^{-1}$, where the \SCCs\ $\tmpcurve_1$ and $\tmpcurve_2$ are oriented, disjoint, and homologous.
\end{thm}
\proof
This theorem was proved for closed surfaces~\cite{Powell} and surfaces with one boundary component~\cite{Johnson}. 
In fact it holds for arbitrary oriented surfaces.

Let $\trldif\in\TrlGrM$. Since $\ubdsurf$ is closed, we have that  
$j(\trldif)$ is generated by diffeomorphisms of types (a) and (b).
Notice that we can choose the corresponding curves so that they belong to $\surf$, 
whence $j(\trldif)$ yields some diffeomorphism $\trldif_1$ of surf such that $\trldif_1^{-1}\circ\trldif \in\ker j$. By Theorem~\ref{th:gen_of_ker_j}, this diffeomorphism is also generated by diffeomorphisms $s_{ik}$ and $r_{ik}$ which evidently are of type (b).
\endproof

\subsection{Proof of \MainTheorem\ for orientable $\surf$ and $\onemanif=\Circle$.}
It suffices to establish the following statement using 
the notations of Lemma~\ref{lm:canon_gen_for_canon_Morse}.
\begin{prop}\label{pr:decomp_diff}
Let $\hdif\in\HMf$ be a diffeomorphism such that
the Morse mappings $\func$ and $\func\circ\hdif:\surf\to\Circle$ are homotopic.
Then $\hdif$ is isotopic to a product of diffeomorphisms of the form
$\orpermdif \circ \cldif \circ \trldif,$
where
\begin{enumerate}
\item
 $\orpermdif$ is generated by $\DifOr$ and those 
  $\bndif_{ij}$ that belong $\HMf$;
\item
 $\cldif$ is generated by Dehn twists along the \SCCs\ of configuration $\CC$ but $\Dtw{\bcurve_{1}}$;
\item
 $\trldif\in\TrlGrM$.
\end{enumerate}
Diffeomorphisms of types {\em (1)-(3)} are $\func$-admissible, whence so is $\hdif$.
\end{prop}
\proof
Evidently that $\hdif$ can be represented as a product $\orpermdif \circ \hdif_1$, where $\hdif_1\in\PHMbd$ and preserves orientation of $\surf$ and $\orpermdif$ is of type (1).
Then by Theorem~\ref{thm:gen_or}, that $\hdif_1$ is generated by the Dehn twists along the curves of configuration $\CC$.

Notice that $\func$ and $\func\circ\hdif_1$ are homotopic.
This condition will allow us to remove $\Dtw{\bcurve_1}$ from the generators of $\hdif_1$ and replace this twist by diffeomorphisms of type (3).

\begin{lem}\label{lm:constr_cldif}
Let $\hdif_1$ be a diffeomorphism of $\surf$
generated by the Dehn twists along the \SCCs\ of configuration $\CC$ and such that $\func$ and $\func\circ\hdif_1$ are homotopic.
Then there exists an $\func$-admissible diffeomorphism $\cldif$ generated by the Dehn twists along the \SCCs\ of configuration $\CC$ except for $\Dtw{\bcurve_1}$
such that the diffeomorphism $\trldif = \cldif^{-1}\circ\hdif_1$ belongs to $\TrlGrM$.
\end{lem}

Hence it remains to establish that every diffeomorphism $\trldif\in\TrlGrM$ is $\func$-admissible.
By Theorem~\ref{th:gen_Torel_grp} it suffices to prove this for diffeomorphism of types (a) and (b).
\begin{thm}\label{th:admis_a_b}
Let $\func:\surf\to\Circle$ be a Morse mapping.

{\rm (i)} Let $\tmpcurve\subset\surf$ be an \SCC\ and $\Dtwcurv$ be a Dehn twist along $\tmpcurve$. 
Then $\Dtwcurv$ is $\func$-admissible iff the restriction $\func|_{\tmpcurve}$ is null-homotopic.
If $\tmpcurve$ separates $\surf$, then $\func|_{\tmpcurve}$ is null-homotopic, whence every diffeomorphism of type {\rm (a)} is $\func$-admissible.

{\rm (ii)}
Every diffeomorphism of type {\rm (b)} is $\func$-admissible.
\end{thm}

Thus in order to complete our proposition, and therefore \MainTheorem, it remains to prove Theorem~\ref{th:admis_a_b} (sections~\ref{sect:prf_th_admis_i} and~\ref{sect:prf_th_admis_ii}) and
Lemma~\ref{lm:constr_cldif} (section~\ref{sect:prf_lm_reduce_torel}). 


\section{Symplectic group}\label{sect:symplect_grp}
For the proof of Lemma~\ref{lm:constr_cldif} we need a description of generators of stabilizers in symplectic group $\SpZg$.
The representation of the group $\SpZg$ is given in~\cite{Birman:SiegelModGr}.
We will also use the ideas from~\cite{OMeara}. 

Let $\modZg$ be a free $2\genus$-module with basis
\chEq\begin{equation}\label{equ:SL_canon_base}
 \acurve_1, \ldots, \acurve_g, \bcurve_1, \ldots, \bcurve_g,
\end{equation}
$I$ be the unity $\genus\times\genus$-matrix,
and $e_{ij}$ be a $\genus\times\genus$-matrix, whose
$(i,j)$-element (the intersection of $i$-th row and $j$-th column)
is equal to $1$ and all other entries are zeros.

Let also $\intform{\cdot}{\cdot}$ be a skew-symmetric $2$-form whose matrix
in the basis~\eqref{equ:SL_canon_base} is the following:
\chEq\begin{equation}\label{equ:SL_matr_std_form}
\left(
\begin{array}{c|c}
0 & I  \\ \hline -I & 0  \\
\end{array}
\right).
\end{equation}
Thus $\intform{\acurve_i}{\bcurve_i} = 1$ and
$\intform{\acurve_i}{\acurve_j}=\intform{\bcurve_i}{\bcurve_j}=\intform{\acurve_i}{\bcurve_j}=0$
for $i,j=1,\ldots,\genus$.
The group of all linear isomorphisms of $\modZg$ preserving $\intform{\cdot}{\cdot}$
is denoted by $\SpZg$ and is called {\em symplectic.}

\subsection{Transvections}
For every $\tmpcurve\in\modZg$ the following automorphism $\Dtw{\tmpcurve}$ of $\modZg$ defined by the formula:
\chEq\begin{equation}\label{equ:transv_def}
\Dtw{\tmpcurve}(x) = \intform{\tmpcurve}{x}\cdot \tmpcurve + x, \qquad \forall x\in\modZg.
\end{equation}
is called the {\em transvection} along $\tmpcurve$.
It is easy to see that $\Dtw{\tmpcurve}\in\SpZg$ and
$$
\Dtw{\tmpcurve}^{-1}(x) = -\intform{\tmpcurve}{x}\cdot \tmpcurve + x, \qquad \forall x\in\modZg.
$$
Define the following elements of $\SpZg$:
\chEq\begin{equation}\label{equ:mu_nu_eta}
\begin{array}{c}
\aastddif{i}{j} = \Dtw{\acurve_i} \circ \Dtw{\acurve_j} \circ  \Dtw{\acurve_i + \acurve_j}^{-1},
\qquad
\bbstddif{i}{j} = \Dtw{\bcurve_i} \circ \Dtw{\bcurve_j} \circ  \Dtw{\bcurve_i + \bcurve_j}^{-1}, \\ [2mm]
\abstddif{i}{j} = \Dtw{\acurve_i} \circ \Dtw{\bcurve_j} \circ  \Dtw{\acurve_i + \bcurve_j}^{-1}.
\end{array}
\end{equation}

\begin{lem}\label{lm:Transvect_formulas}
The following formulas hold true
for $i\not=j = 1\ldots\genus$:
\begin{gather*}
\Dtw{\acurve_i} =
\mybrack{
\begin{array}{c|c}
I & \Ematr{i}{i}  \\ \hline 0 & I  \\
\end{array}
}, \enspace
\Dtw{\bcurve_i}=
\mybrack{
\begin{array}{c|c}
I & 0 \\ \hline -\Ematr{i}{i}  & I  \\
\end{array}
}, \enspace
\Dtw{\acurve_i + \bcurve_j} =
\mybrack{
\begin{array}{c|c}
I-\Ematr{i}{j} & \Ematr{i}{i} \\ \hline
-\Ematr{j}{j}  & I + \Ematr{j}{i} \\
\end{array}
}.
\\
\Dtw{\acurve_i + \acurve_j} =
\mybrack{
\begin{array}{c|c}
I & \Ematr{i}{i} + \Ematr{j}{j} + \Ematr{i}{j} + \Ematr{j}{i}\\ \hline  0 & I  \\
\end{array}
},
\\
\Dtw{\bcurve_i + \bcurve_j} =
\mybrack{
\begin{array}{c|c}
I & 0 \\ \hline
-\Ematr{i}{i} - \Ematr{j}{j} - \Ematr{i}{j} - \Ematr{j}{i} & I  \\
\end{array}
},
\\
\aastddif{i}{j} =
\mybrack{
\begin{array}{c|c}
I & - \Ematr{i}{j} - \Ematr{j}{i} \\ \hline  0 & I  \\
\end{array}
}, \enspace
\bbstddif{i}{j} =
\mybrack{
\begin{array}{c|c}
I & 0 \\ \hline \Ematr{i}{j} + \Ematr{j}{i} & I  \\
\end{array}
} \\ 
\abstddif{i}{j} =
\mybrack{
\begin{array}{c|c}
I + \Ematr{i}{j} & 0 \\ \hline  0 & I - \Ematr{j}{i} \\
\end{array}
}.
\end{gather*}
Moreover, the matrices
$\Dtw{\acurve_i}$, $\Dtw{\bcurve_i}$, $\aastddif{i}{j}$, $\bbstddif{i}{j}$, and $\abstddif{i}{j}$,
$(i\not=j=1\ldots\genus)$ generate $\SpZg$.
\end{lem}
\proof
The lemma can be established by direct calculations.
The fact that these matrices generate $\SpZg$ can be easily deduced from~\cite[Ch.~2, \S2.2.]{OMeara} or~\cite{Birman:SiegelModGr}.
\endproof

For each $x\in\modZg$ denote by $\TransAnn(x)$
the subgroup in $\SpZg$ generated by transvections along elements of $\modZg$
that are $\intform{\cdot}{\cdot}$-orthogonal to $x$, i.e.
\chEq\begin{equation}\label{equ:TransAnn}
 \TransAnn(x) = \langle \Dtw{\tmpcurve} \ | \ \tmpcurve\in\modZg, \intform{\tmpcurve}{x}=0  \rangle.
\end{equation}
Let also $\Stab(x)$ be the stabilizer of $x$ in $\SpZg$, i.e.
$$\Stab(x) = \{\hdif\in\SpZg \ | \ \hdif(x)= x\}.$$
It easily follows from~\eqref{equ:transv_def} that $\TransAnn(x)\subset\Stab(x)$.

\begin{prop}\label{pr:TransAnn_a_1}
$\TransAnn(\acurve_1) = \Stab(\acurve_1)$. Moreover, this group is generated by the following matrices:
\begin{multline}\label{equ:TransAnn_a1_gen}
\Dtw{\acurve_i}, \ \Dtw{\bcurve_i}, \ \aastddif{i}{j}, \ \bbstddif{i}{j}, \ \abstddif{i}{j},
\\ \text{except for} \ \Dtw{\bcurve_1}, \ \bbstddif{1}{i}=\bbstddif{i}{1} \  \text{and\/}
\ \abstddif{i}{1}, \
(i\not=j=1,\ldots,\genus).
\end{multline}
\end{prop}
\proof
Evidently, the matrices~\eqref{equ:TransAnn_a1_gen} belong to $\TransAnn(\acurve_1)$.
Let $\hdif\in\Stab(\acurve_1)$.
We will show that $\hdif$ is generated by~\eqref{equ:TransAnn_a1_gen}.
The proof consists of two steps.

{\em Step 1.}
We will find an element $\hdif_1\in\SpZg$ such that $\hdif\cdot\hdif^{-1}_1$
is generated by~\eqref{equ:TransAnn_a1_gen} and $\hdif_1(\bcurve_1)=\bcurve_1$.
Let
$$
\hdif(\bcurve_1) =
\acoef_1\,\acurve_1 + \bcoef_1\,\bcurve_1 +
\acoef_2\,\acurve_2 + \bcoef_2\,\bcurve_2 + \ldots,
$$
for some $\acoef_i, \bcoef_i\in\ZZZ, (i=1\ldots,\genus)$.
Since $\hdif$ preserves the form $\intform{\cdot}{\cdot}$
and $\hdif(\acurve_1)=\acurve_1$, we get
$$ \bcoef_1=
\intform{\acurve_1}{\hdif(\bcurve_1)} =
\intform{\hdif(\acurve_1)}{\hdif(\bcurve_1)}=
\intform{\acurve_1}{\bcurve_1}=1.
$$

Consider now the effect of action of $\aastddif{1}{j}$ and $\abstddif{1}{j}$ on
$\hdif(\bcurve_1)$, $j=2\ldots\genus$.
Let $t\in\ZZZ$. Then it is easy to verify that for $j>1$ we have:
\begin{gather*}
(\aastddif{1}{j})^{t} \circ \hdif (\bcurve_1) =
 (\acoef_1 - t\bcoef_j)\,\acurve_1 + \bcurve_1 + \ldots + (\acoef_j-t)\,\acurve_j +
                     \bcoef_j\,\bcurve_j + \ldots,
\\
(\abstddif{1}{j})^{t} \circ \hdif (\bcurve_1) =
 (\acoef_1+t\acoef_j)\,\acurve_1 + \bcurve_1 + \ldots + \acoef_j\,\acurve_j + (\bcoef_j-t)\,\bcurve_j + \ldots,
\end{gather*}
where the coefficients at other basis elements are not changed.

Define now $\hdif_1\in\modZg$ by the formula:
$$ \hdif_1 =
(\Dtw{\acurve_1})^{-\acoef'}
\cdot \prod_{j=2}^{\genus} (\abstddif{1}{j})^{\bcoef_j}  \cdot
 \prod_{i=2}^{\genus} (\aastddif{1}{j})^{\acoef_j} \cdot \hdif,
$$
where
$$\acoef' = {\acoef_1 - \sum\limits_{j=2}^{g} \acoef_j \bcoef_j}.$$
We claim that $\hdif_1(\bcurve_1)=\bcurve_1$.

Indeed, the product of $\aastddif{1}{j}$ reduces the coefficients at $\acurve_j$
and the product of $\abstddif{1}{j}$ reduces the coefficients at $\bcurve_j$ for every $j=2\ldots\genus$.
This also makes the coefficient at $\acurve_1$ equal to $\acoef'$.
Since
$$
\Dtw{\acurve_1}(\acurve_1)=\acurve_1  \qquad \text{and} \qquad
(\Dtw{\acurve_1})^{t}(\bcurve_1) = (\acoef_1 + t)\,\acurve_1 + \bcurve_1,
$$
we obtain that the multiple $(\Dtw{\acurve_1})^{-\acoef'}$ reduces this coefficient.

{\em Step 2.}
Consider the following submodules of $\modZg$:
$$
P= \langle \acurve_1, \bcurve_1 \rangle \qquad \text{and} \qquad
Q = \langle \acurve_i, \bcurve_i \ | \ {i=2\ldots\genus} \rangle.
$$
They are orthogonal with respect to the form $\intform{\cdot}{\cdot}$ and $\hdif_1|_{P} = \id$.
Since $\hdif_1$ preserves $\intform{\cdot}{\cdot}$, it follows that $\hdif_1(Q)=Q$.
Thus $\hdif_1$ can be regarded as an element of
the group $\Sp{2\genus-2}{\ZZZ} \subset \Sp{2\genus}{\ZZZ}$
consisting of isomorphisms that are identity on $P$.

By Lemma~\ref{lm:Transvect_formulas} the group
$\Sp{2\genus-2}{\ZZZ}$ is generated by matrices~\eqref{equ:TransAnn_a1_gen} for $i\not=j=2\ldots\genus$.
In particular, they generate $\hdif_1$.
\endproof


\section{Minimal Morse maps.}\label{sect:minimal_Morse_maps}
For the proof of Theorem~\ref{th:admis_a_b} we need the notion of minimal Morse mapping.
Let $\surf$ be a compact surface orientable or not. 
We say that a Morse map $\func:\surf\to\onemanif$ is {\em minimal\/}
if the number $\crpntf{0}+\crpntf{1}+\crpntf{2}$ of critical points of $\func$
is minimal among all possible Morse maps $\surf\to\onemanif$ having the same sets of positive and negative boundary components as $\func$.
Let $\bndposf$ and $\bndnegf$ be the number of $\func$-positive and $\func$-negative boundary components of $\surf$.
The following lemma is easy to prove:
\begin{lem}\label{lm:min_func_crit}
A Morse mapping $\func:\surf\to\onemanif$ is minimal iff for every connected component $\comp$ of $\surf$ the restriction $\func|_{\comp}$ is minimal.
A Morse function $\func:\surf\to\Rline$ on a {\em connected\/} surface $\surf$
is minimal if and only if the following two relations hold true
\chEq\begin{equation}\label{equ:min_func_crit}
\crpnt{0}{\func}=\left\{
\begin{array}{ll}
1, & \mathrm{if} \ \bndnegf=0, \\
0, & \mathrm{if} \ \bndnegf>0,
\end{array}
\right.
\quad
\crpnt{2}{\func}=\left\{
\begin{array}{ll}
1, & \mathrm{if} \ \bndposf=0, \\
0, & \mathrm{if} \ \bndposf>0.
\end{array}
\right.
\end{equation}
Let $\func:\surf\to\Circle$ be a Morse mapping which is not null-homotopic.
Then $\func$ is minimal iff \ $\crpntf{0}=\crpntf{2}=0$. \qed
\end{lem}

We admit now that $\surf$ may be not connected.
Let $\func:\surf\to[0,1]$ be a Morse function such that $\onehalf \in[0,1]$ is its regular value.
Denote
$$\surfp_0=\func^{-1}\left[0,\Onehalf\right], \quad \surfp_1=\func^{-1}\left[\Onehalf,1\right].$$
$$ \btsurf{0}=\func^{-1}(0), \quad \btsurf{1} = \func^{-1}(1), \quad \Invim=\func^{-1}\left(\Onehalf\right)$$

\begin{lem}\label{lm:min_Morse_maps}
Suppose that 
\begin{enumerate}
\item
$\btsurf{0}$, $\btsurf{1}$ and $\Invim$ are nonempty, the union $\btsurf{0}\cup\btsurf{1}$ is included in $\partial\surf$ and intersects every connected component of $\surf$;

\item
the restriction $\func|_{\surfp_i}$ is a {\em minimal\/} Morse function for $i=0,1$;

\item
for every connected component $\comp$ of $\surf$ such that $\comp\cap\Invim\not=\varnothing$ 
we have $\comp \cap \btsurf{i} \not=\varnothing$ for both $i=0,1$.
\end{enumerate}
Then $\func$ is a minimal Morse function on $\surf$.
\end{lem}
\proof
Let $\comp$ be a component of $\surf$.
We will show that $\func|_{\comp}$ is a minimal Morse function.
Denote $\comp_i=\comp\cap\surfp_i$ $(i=0,1)$.

If $\comp\cap\Invim=\varnothing$, then $\comp$ is a connected component of
one of the sets either $\surfp_0$ or $\surfp_1$.
Then the restriction of $\func$ onto $\comp$ is minimal.

Suppose that $\comp\cap\Invim\not=\varnothing$.
Then $\comp\cap\btsurf{i}\not=\varnothing$ for $i=0,1$ by (3).
Evidently, the components of the intersection $\comp\cap\Invim\not=\varnothing$ are negative
for the restriction $\func|_{\comp_{1}}$
and positive for the restriction $\func|_{\comp_{0}}$.
Therefore, by Lemma~\ref{lm:min_func_crit}, we have
\chEq\begin{equation}\label{equ:crpnt_count_1}
 \crpnt{2}{\func|_{\comp_{0}}} = \crpnt{0}{\func|_{\comp_{1}}} = 0.
\end{equation}

Similarly, the intersection $\comp\cap\btsurf{0}$ (resp. $\comp\cap\btsurf{1}$)
consists of some negative (resp. positive) components of $\func|_{\comp}$
and $\func|_{\comp_0}$ (resp. $\func|_{\comp_1}$).
Then from Lemma~\ref{lm:min_func_crit}, we also get
$$
 \crpnt{0}{\func|_{\comp_{0}}} = \crpnt{2}{\func|_{\comp_{1}}} = 0.
$$
Combining this with~\eqref{equ:crpnt_count_1}, we obtain
$$
\crpnt{i}{\func|_{\comp}} = \crpnt{i}{\func|_{\comp_{0}}} + \crpnt{i}{\func|_{\comp_{1}}} = 0, \quad i=0,2.
$$
Whence by Lemma~\ref{lm:min_func_crit} $\func|_{\comp}$ is minimal.
\endproof


\section{Minimization of intersections with a level-set}\label{sect:minimizat_int_levsets}
Let $\surf$ be a compact surface (orientable or not),
$\func:\surf\to\Circle$ be a Morse mapping, and $\tmpcurve_1,\ldots,\tmpcurve_m\subset\surf$ be disjoint \SCCs.

\begin{lem}\label{lm:curv_min_int}
$\func$ is \sgmhom-homotopic to a Morse mapping $\gfunc$ such for some level-set $\lset$ of $\gfunc$ and for every $i=1,\ldots,m$ the curve $\tmpcurve_i$ does not pass through the critical points of $\gfunc$ and 
\begin{enumerate}
\item[(i)]
if the restriction $\func|_{\tmpcurve_i}$ is not null-homotopic, then $\tmpcurve_i$ transversely intersects every level-set of $\gfunc$;

\item[(ii)]
otherwise $\tmpcurve_i \cap \lset = \varnothing$.
\end{enumerate}
\end{lem}
\proof
Let $\cval\in\Circle$ be a regular value of $\func$.
Set 
$$\gunion = \mathop\cup\limits_{i=1}^{m}\tmpcurve_i, \qquad
\indnum=\concomp{\func^{-1}(\cval)  \cap \gunion},
\qquad \text{and} \qquad
\sumdeg=\sum_{i=1}^{m} |\deg\func|_{\tmpcurve_i}|.$$
Then $\concomp{\func^{-1}(\cval)  \cap \tmpcurve_i} \geq \deg\func|_{\tmpcurve_i}$
for $i=0,1$, whence $\indnum \geq \sumdeg$.
Moreover, $\indnum = \sumdeg$ if and only if
$\concomp{\func^{-1}(\cval)  \cap \tmpcurve_i} = \deg\func|_{\tmpcurve_i}$.

\begin{claim}\label{clm:res_min_intersect}
Suppose that $\indnum > \sumdeg$.
Then $\func$ is \sgmhom-homotopic to a Morse map $\func_1$ such that
$\concomp{\func_1^{-1}(\cval_1) \cap \gunion} < \indnum$
for some regular value $\cval_1$ of $\func_1$.
\end{claim}
\proof
We will exploit the notations and the construction of Section~\ref{subsect:CuttingAlongLevelSet}.
Cutting $\surf$ along $\func^{-1}(\cval)$ we obtain the surface
$\tsurf$ and the Morse function $\tfunc:\tsurf\to[0,1]$.
Let also $\psurf:\tsurf\to\surf$ be the factor-map,
$\btsurf{i}=\tfunc^{-1}(i)$ for $i=0,1$,
and $\btsurf{}=\btsurf{0}\cup\btsurf{1}=\psurf^{-1}(\func^{-1}(\cval))$.

Let $\Arcs = \psurf^{-1}(\gunion)$ and $\arc_1,\ldots,\arc_{\narc}$ be the connected components of $\Arcs$.
Then the intersection $\arc_j\cap\btsurf{}$ is either empty (whence $\arc_j$ is an \SCC)
or consists of two points (whence $\arc_j$ is a simple arc with ends in $\btsurf{}$).
Let us divide $\Arcs$ into four groups $\arccompE$, $\arccompAA$, $\arccompAB$, $\arccompBB$
consisting of arcs
that respectively do not intersect $\btsurf{}$,
intersect only $\btsurf{0}$,
intersect both sets $\btsurf{0}$ and $\btsurf{1}$,
and intersect only $\btsurf{1}$.
Thus $\Arcs = \arccompE \cup \arccompAA \cup \arccompAB \cup \arccompBB$.
Notice that $\concomp{\Arcs\cap\btsurf{0}} = \concomp{\Arcs\cap\btsurf{1}} = \indnum$,
$\concomp{\arccompAA} = \concomp{\arccompBB}$,
and the sets $\arccompAA$ and $\arccompBB$ are non-empty if and only if
$\indnum > \sumdeg$.

Let $\compAB\subset\tsurf$ be a union of those connected components of $\tsurf$
which intersect both sets $\btsurf{0}$ and $\btsurf{1}$.
Consider the set
$$
  G=\compAB \cap (\btsurf{0}\cup\arccompAA).
$$
By definition, $G\cap (\arccompE\cup\arccompBB)=\varnothing$.
Then there exists a regular neighborhood $\aNbh$ of $G$ which does not intersect $\arccompE\cup\arccompBB$ and such that the boundary $\Invim = \partial\aNbh$ transversely intersects every component of $\arccompAB$ at a unique point.
Hence, $\Invim\cap\Arcs = \Invim \cap \arccompAB$.
Evidently, $\Invim$ separates $\tsurf$ between $\btsurf{0}$ and $\btsurf{1}$.
Moreover, $\concomp{\Invim \cap \arccompAB} < \indnum$.

We will now construct
a Morse function $\tgfunc:\tsurf\to[0,1]$ which
coincides with $\tfunc$ in some neighborhood of $\btsurf{}\cup\partial\tsurf$,
has critical type of $\tfunc$, and such that $\tgfunc^{-1}(\onehalf)=\Invim$.

Let $\tgfunc_0:\Nbh_0\to[0,\onehalf]$ and
$\tgfunc_1:\Nbh_1\to[\onehalf,1]$ be two minimal Morse functions such that
$$
\tgfunc_0^{-1}(0)=\btsurf{0},
\qquad
\tgfunc_0^{-1}\left(\Onehalf\right)=\tgfunc_1^{-1}\left(\Onehalf\right) = \Invim,
\qquad
\tgfunc_1^{-1}(1)=\btsurf{1},
$$
and the Morse function $\tgfunc:\tsurf\to[0,1]$
defined by $\tgfunc|_{\Nbh_i}=\tgfunc_i$, $(i=0,1)$ is $C^{\infty}$,
has the same sets of positive and negative components as $\tfunc$,
and coincide with $\tfunc$ in some neighborhood of $\btsurf{}\cup\partial\tsurf$.

We claim that $\tgfunc$ is minimal.
Indeed, let $\comp$ be a component of $\tsurf$ such that $\comp\cap\Invim\not=\varnothing$.
Since $\Invim=\partial\aNbh\subset \compAB$, we obtain that $\comp\subset\compAB$.
Denote $\comp_i = \comp \cap \Nbh_i$,
then $\comp\cap\btsurf{i}=\comp_{i}\cap\btsurf{i}\not=\varnothing$,
by the definition of $\compAB$.
It follows from Lemma~\ref{lm:min_Morse_maps} that $\tgfunc$ is minimal.

Adding critical points to $\tgfunc$ outside of $\btsurf{}\cup\Invim$
we can change its critical type to the critical type of $\tfunc$.
Let us denote this new function by $\tfunc_1$.
Then $\tfunc_1$ satisfies the statement of our claim.

Indeed, denote $\cval_1=\pcirc(\onehalf)$.
By the case $\onemanif=\Rline$ of \MainTheorem\ we obtain that $\tfunc\sgmh\tfunc_1$ with respect to some neighborhood of $\btsurf{}\cup\partial\tsurf$.
This \sgmhom-homotopy induces a \sgmhom-homotopy (with respect to $\func^{-1}(\cval)\cup\partial\surf$)
of $\func$ to a Morse mapping $\func_1$ such that
$\concomp{\func_1^{-1}(\cval_1)  \cap \gunion} <\indnum$.
\endproof

We now proceed with the proof of Lemma~\ref{lm:curv_min_int}.  
By Claim~\ref{clm:res_min_intersect} we can assume that $\indnum = \sumdeg$.
As noted above this is equivalent to the statement:
$\concomp{\func^{-1}(\cval)\cap\tmpcurve_i}=\deg\func|_{\tmpcurve_i}$.
In particular, if the restriction $\func|_{\tmpcurve_i}$ is null-homotopic,
then $\concomp{\func^{-1}(\cval)\cap\tmpcurve_i}=0$, i.e.
$\tmpcurve_i\cap\func^{-1}(\cval)=\varnothing$, whence (ii) holds true.

Let us assume that $\arc_i$ is given by an embedding $\arc_i:[0,1]\to\tsurf$
so that $\arc_i\cap\arc_j=\varnothing$ for $j\not=i$.
To establish (i) we prove that following claim:
\begin{claim}\label{clm:curv_monotone}
Suppose that $\arc_i(0)\in\btsurf{0}$, $\arc_i(1)\in\btsurf{1}$,
and the intersection $\arc_i\cap\btsurf{}$ is transversal for each $i=1,\ldots,\narc$.
Then $\tfunc$ is \sgmhom-homotopic to a Morse function $\tgfunc$ such
that $\arc_i$ is transversal to level-sets of $\tgfunc$.
\end{claim}
It follows that a \sgmhom-homotopy of this claim
yields a \sgmhom-homotopy $\func\sgmh\gfunc$ with respect $\func^{-1}(\cval)$
such that every $\tmpcurve_i$ is transversal to level-sets of $\gfunc$.
This will complete Lemma~\ref{lm:curv_min_int}.
\proof[Proof of Claim~\ref{clm:curv_monotone}]
We will construct a Morse function $\tfunc_1$ and a gra\-dient-like vector field $\Fld$ for $\tfunc_1$ such that for every $i=1,\ldots,m$ the arc $\arc_i$ is a trajectory of $\Fld$.
Then adding or canceling the proper number of pairs of critical points of $\tfunc_1$ outside $\cup_i\arc_i$ we obtain a Morse function $\tgfunc$ having the critical type of $\tfunc$ and such that $\Fld$ is a gradient-like for $\tgfunc$.

For every $i=1,\ldots,m$ let
$\phi_i:[0,1]\times[-1,1]\to \manif$ be a smooth  embedding
such that the image $\Nbh_i=\IM\phi_i$ is a neighborhood of $\arc_i$,
$\phi_i(t,0) = \arc_i(t)$ for $t\in[0,1]$,
$\phi^{-1}(\btsurf{s}) = \{s\}\times[-1,1]$ for $s=0,1$.
Since $\arc_i$ are mutually disjoint, we can assume that so are $\Nbh_i$.
Denote $\Nbh = \cup_{i=1}^{m}\Nbh_i$ and define
a function $\tgfunc:\Nbh\to[0,1]$ by the formula
$\tgfunc(x)=p_2 \circ\phi_i^{-1}(x)$ for $x\in\Nbh_i$,
where $p_2:[0,1]\times[-1,1]\to[-1,1]$ is the natural projection.

Slightly changing $\tgfunc$ outside some neighborhood of $\cup_{i}\arc_i$ we can extend $\tgfunc$ over all of $\tsurf$.
Moreover, this extension may be assumed Morse whose positive and negative boundary components coincide with ones of $\tfunc$ though the number of critical points of $\tgfunc$ and $\tfunc$ may be different. 
Now we show how to change the critical type $\Crtype{\tgfunc}$ of $\tgfunc$ to $\Crtype{\tfunc}$ by adding or canceling pairs of critical points outside of $\cup_{i}\arc_i$.

Recall that a vector field $\Fld$ on a manifold $\tsurf$ is {\em gradient-like\/}
for a function $\tfunc:\tsurf\to\Rline$ if $d\tfunc(\Fld)(x)>0$
at each regular point $x$ of $\tfunc$.

Let $\gFld$ be any gradient-like vector field for the function $\tgfunc$ on $\tsurf$
and $\tpFld$ be the gradient vector field
for the function $p_2$ on $[0,1]\times[-1,1]$, i.e. $\tpFld(s,t)=(0,1)$.
Using $\phi_i$ we transfer $\tpFld$ to $\Nbh_i$.
This gives us a vector field $\pFld$ on $\Nbh$
such that $\arc_i$ is a trajectory of $\pFld$ for $i=1,\ldots,m$.

Finally, we glue $\gFld$ and $\pFld$.
Let $\Nbh'$ be a neighborhood of $\cup_{i}\arc_i$ such that
$\overline{\Nbh'} \subset \Nbh$ and let $\aNbh=\tsurf\setminus\overline{\Nbh'}$.
Then $\Nbh\cup\aNbh=\tsurf$.

Let $\mu_1,\mu_2:\tsurf\to[0,1]$ be a partition of unity
corresponding to the open covering $\{\Nbh, \aNbh\}$ of $\tsurf$, i.e.
$\supp\mu_1\subset\Nbh$,
$\supp\mu_2\subset\aNbh$,
and $\mu_1+\mu_2\equiv 1$.
Define a vector field $\Fld$ on $\tsurf$ by the formula
$$ \Fld(x) = \mu_1(x)\cdot \pFld(x) +  \mu_2(x)\cdot \gFld(x), \qquad x\in\tsurf.$$
Evidently, $\Fld$ is gradient-like for $\tgfunc$
and coincides with $\pFld$ near $\cup_{i}\arc_i$.
In particular, every $\arc_i$ is a trajectory of $\Fld$,
whence $\arc_i$ transversely intersects level-sets of $\tgfunc$.

It remains to show that $\tgfunc$ can be changed outside $\cup_{i}\arc_i$
to have critical type of $\tfunc$.
First we show how to make $\tgfunc$ a minimal Morse function.

Suppose that $\tgfunc$ has a critical point $z_0$ either of index $0$ or $2$.
Since the sets of positive and negative boundary components of $\tgfunc$ are non-empty,
there exists a critical points $z_1$ of index $1$
and a trajectory $\omega$ of $\Fld$ with ends at $z_0$ and $z_1$.
This trajectory does not intersect $\cup\arc_i$.
Hence $\tgfunc$ can be changed in some neighborhood $N$ of $\omega$
to have no critical points in $N$ (see~\cite{Morse_Huebsch, Milnor_h_cob}).
Thus the number of critical points is reduced.
By the similar procedure we can add pairs of critical points outside of $\cup_i \arc_i$.
Therefore we can change the critical type $\Crtype{\tgfunc}$ of $\tgfunc$ to $\Crtype{\tfunc}$ leaving $\arc_i$ transversal to level-sets of $\tgfunc$.
\endproof


\section{Proof of (i) of Theorem~\ref{th:admis_a_b}}\label{sect:prf_th_admis_i}
Let $\tmpcurve\subset\surf$ be a simple closed curve
and $\Dtwcurv$ be a Dehn twist along $\tmpcurve$.

{\em Necessity.}
Suppose that $\Dtwcurv$ is $\func$-admissible. 
Then $\func$ and $\func\circ\Dtwcurv$ are homotopic.
We should show that $\deg\func|_{\tmpcurve}=0$.
We can assume that there is a regular value $\cval$ of $\func$ such that $\acurve=\func^{-1}(\cval)$ is an \SCC.
Denote $\acurve' = \Dtwcurv(\acurve)$.

Since $\func$ and $\func\circ\Dtwcurv$ are homotopic, we obtain from the last paragraph of Section~\ref{sect:or_level_sets} that $[\acurve']=[\acurve]$ in $H_1(\surf,\partial\surf)$, i.e.\! $\Dtwcurv$ fixes $[\acurve]$.
Then by Eq.~\eqref{equ:transv_def} for the action of Dehn twist in $H_1(\surf,\partial\surf)$ we get
$$ [\acurve] = \Dtwcurv([\acurve]) = \intform{[\tmpcurve]}{[\acurve]}\cdot [\tmpcurve] + [\acurve] = \deg\func|_{\tmpcurve}\cdot[\tmpcurve] + [\acurve],$$
whence $\deg\func|_{\tmpcurve}=0$.

{\em Sufficiency.}
Suppose that $\func|_{\tmpcurve}$ is null-homotopic.
By Lemma~\ref{lm:curv_min_int}, $\func$ is \sgmhom-homotopic to a Morse mapping $\gfunc$ such that
$\gfunc^{-1}(\cval)\cap\tmpcurve=\varnothing$ for some regular value $\cval$ of $\gfunc$.
We apply now the construction of Section~\ref{subsect:CuttingAlongLevelSet}.
Cutting $\surf$ along $\gfunc^{-1}(\cval)$ we obtain a surface $\tsurf=\tsurf\frval{\gfunc}{\cval}$,
a Morse function $\tgfunc:\tsurf\to[0,1]$, and an \SCC\  $\tcurve\subset\tsurf$
corresponding to $\tmpcurve$.
By the case $\onemanif=\Rline$ of \MainTheorem, $\Dtw{\tcurve}$ is $\tgfunc$-admissible.
Then $\Dtwcurv$ is $\gfunc$-admissible and therefore is $\func$-admissible.
\endproof


\section{Proof of (ii) of Theorem~\ref{th:admis_a_b}}\label{sect:prf_th_admis_ii}
Let $\func:\surf\to\Circle$ be a Morse mapping,
$\tmpcurve_1$, $\tmpcurve_2$ be disjoint oriented homologous simple closed curves in $\surf$, and $\Dtwist=\Dtw{\tmpcurve_1}\circ\Dtw{\tmpcurve_2}^{-1}$ be the product of Dehn twists along these curves.
We must prove that $\Dtwist$ is $\func$-admissible.

Since these curves are homologous, it follows that the restrictions of $\func$ to them are homotopic. 
If these restrictions are null-homotopic, then by the case (i) of this theorem $\Dtwist$ is $\func$-admissible.
Therefore we will assume that $\func|_{\tmpcurve_1}\not\sim0$.

By Lemma~\ref{lm:curv_min_int} we can also assume that $\tmpcurve_i$ transversely intersects each level-set of $\gfunc$.
Then the statement (ii) of Theorem~\ref{th:admis_a_b} is a direct corollary of the following lemma:

\begin{lem}\label{lm:sgmhom_gfunc_gfunc_tw}
Let $\func:\surf\to\Circle$ be a Morse mapping,
$\tmpcurve_1$, $\tmpcurve_2$ be two disjoint homologous \SCCs\ in $\surf$,
and $\Dtwist = \Dtw{\tmpcurve_1} \circ \Dtw{\tmpcurve_2}^{-1}$.
Suppose that both of $\tmpcurve_i$ transversely intersect every level-set of $\func$.
Then $\func\sgmh\func\circ\Dtwist$.
\end{lem}
\proof
Let $\comp \subset \surf$ be the closure of one of the
connected components of $\surf\setminus(\tmpcurve_1 \cup \tmpcurve_2)$
bounded by the curves $\tmpcurve_1$ and $\tmpcurve_2$.
Since $\tmpcurve_{\indk}$ $(\indk=1,2)$ transversely intersects level-sets of $\gfunc$,
there exists an embedding $\phi_{\indk}$ of $\Circle\times[-2,2]$ onto
some neighborhood $\collar_{\indk}$ of $\tmpcurve_{\indk}$ such that
\chEq\begin{equation}\label{equ:cond_on_phi_ind}
\phi_{\indk}(\Circle\times\{0\})=\tmpcurve_{\indk}, \qquad
\phi_{\indk}(\Circle\times[0,2])\subset\comp,
\end{equation}
and the following diagram is commutative:
\chEq\begin{equation}\label{equ:CD_g_phi_zd}
\begin{CD}
\Circle\times[-2,2] @>{\phi_{\indk}}>> \collar_{\indk} & \ \subset \ & \surf \\
@V{p_1}VV @VV{\gfunc}V \\
\Circle @>{\cov}>> \Circle \\
\end{CD}
\end{equation}
Here $p_1$ is a projection onto the first coordinate
and $\cov$ is a covering mapping of degree $d = \deg\func|_{\tmpcurve_1} = \deg\func|_{\tmpcurve_2}$
defined by the formula $\cov(z)=z^d$.
Thus
\chEq\begin{equation}\label{equ:g_phi_zd}
\gfunc\circ \phi_{\indk}(z,t) = z^d.
\end{equation}
We can also assume that $\collar_1\cap\collar_2=\varnothing$.
To simplify notations for each pair $a,b\in[-2,2]$ we denote
$$ \collar_{\indk}^{[a,b]} = \phi_{\indk}(\Circle\times[a,b]).$$

Let $\mu:[-2,2]\to[0,1]$ be a $C^{\infty}$ function such that $\mu[-2,-1] = 0$ and $\mu[1,2]=1$.
Then the Dehn twist $\Dtw{\tmpcurve_{\indk}}$ along $\tmpcurve_{\indk}$ can be defined so that $\Dtwist=\Dtw{\tmpcurve_{1}}\circ\Dtw{\tmpcurve_{2}}^{-1} $ will have the form:
\chEq\begin{equation}\label{equ:pair_Dehn_twist}
\Dtwist(z,t) = \left\{
\begin{array}{cl}
x, & x\in\surf\setminus(\collar_1\cup\collar_2) \\
(z e^{2\pi i \mu(s)},s), &  x = \phi_k(z,\sVar)\in\collar_k, k=1,2. \\
\end{array}
\right.
\end{equation}

Now a \sgmhom-homotopy $\Gsgm:\surf\times[0,1]\to\Circle$
between $\gfunc$ and $\gfunc\circ\Dtwist$ can be defined by the formula:
$$
\Gsgm(x,t)=\left\{
\begin{array}{ll}
\gfunc(x)\,e^{2\pi i d t}, & x \in \comp \setminus(\collar_1^{[0,1]} \cup \collar_2^{[0,1]}),\\
\gfunc \circ \phi_k(z\,  e^{2\pi i \mu(\sVar) \cdot t},\sVar), &
  x = \phi_k(z,\sVar)\in\collar_k, k=1,2. \\
\gfunc(x), & x \in \surf\setminus(\comp \cup \collar_1^{[-1,0]} \cup \collar_2^{[-1,0]}).
\end{array}
\right.
$$
\begin{rem}
A geometrical meaning of this formula is that the mapping $\Gsgm$ ``moves'' $d$ times the part $\comp$ between the curves $\tmpcurve_1$ and $\tmpcurve_2$  ``around $\Circle$'' leaving the complement $\surf\setminus\comp$ fixed.
\end{rem}

Let us verify, that $\Gsgm$ is in fact a \sgmhom-homotopy
connecting $\gfunc$ with $\gfunc\circ\Dtwist$.
\proof
It is clear that $\Gsgm_0=\gfunc$.
Moreover, it follows from~\eqref{equ:cond_on_phi_ind} and~\eqref{equ:CD_g_phi_zd}
that $\phi_1$ preserves orientation of $\Circle\times[-2,2]$ while $\phi_2$ reverses it.
Hence by~\eqref{equ:pair_Dehn_twist} we get
$\Gsgm_1 = \gfunc\circ\Dtw{\tmpcurve_1}\circ\Dtw{\tmpcurve_2}^{-1}$.

Evidently, the continuity of $\Gsgm$ will imply its smoothness.
To prove that $\Gsgm$ is continuous we should verify that the second formula
coincides with the first one on $\collar_1^{[1,2]} \cup \collar_2^{[1,2]}$
and with the third one on $\collar_1^{[-2,-1]} \cup \collar_2^{[-2,-1]}$.

Let $x = \phi_k(z,\sVar)\in\collar_k^{[1,2]}$ for $k=1,2$, then $\mu(\sVar)=1$, whence,
using~\eqref{equ:g_phi_zd}, we get
$$
\gfunc \circ \phi_k(z\,  e^{2\pi i \mu(\sVar) \cdot t},\sVar) =
z^d \, e^{2\pi i d t} = \gfunc(x) \, e^{2\pi i d t}.
$$

Let now $x = \phi_k(z,\sVar)\in\collar_k^{[-2,-1]}$ for $k=1,2$, then $\mu(\sVar)=0$, whence
$$\gfunc \circ \phi_i(z\,  e^{2\pi i \mu(\sVar) \cdot t},\sVar)=\gfunc \circ \phi_i(z,\sVar)=\gfunc(x).
$$

Notice that for every point $x\in\manif$ there exists a neighborhood
on which $\Gsgm_t$ differs from $\gfunc$
by a diffeomorphism of either $\Circle$ or $\surf$.
Hence $\Gsgm_t$ is Morse for all $t\in[0,1]$,
i.e. $\Gsgm$ is a \sgmhom-homotopy.
\endproof


\section{Proof of Lemma~\ref{lm:constr_cldif}}\label{sect:prf_lm_reduce_torel}
Suppose that $\hdif\in\PHMbd$ is generated by $\{\Dtw{l}: l\in\CC\}$ and such that the mappings $\func$ and $\func\circ\hdif$ are homotopic.
We have to prove that $\hdif$ is in fact generated by 
$\{\Dtw{l}: l\in\CC\setminus\bcurve_1\}$.

Recall that $H_1(\surf,\partial\surf)$ is a free module generated by homology classes of $\acurve_1, \ldots, \acurve_g, \bcurve_1, \ldots, \bcurve_g$.
Moreover, the matrix of $\intform{\cdot}{\cdot}$ in this basis has the form~\eqref{equ:SL_matr_std_form}. 
Since $\hdif_{\hommap}$ preserves this $\intform{\cdot}{\cdot}$ we may suppose that $\hdif_{\hommap}\in\SpZg$.

Notice that $\hdif_{\hommap}[\acurve_1]=[\acurve_1]$, since $\acurve_1$ is a level-set of $\func$, whence $\hdif_{\hommap}$ belongs to the stabilizer $\Stab([\acurve_1])$ of $\acurve_1$ in $\SpZg$.

Let $\Dtw{\tmpcurve}$ be a Dehn twist along simple closed curve $\tmpcurve$.
Then it acts on $H_1(\surf,\partial\surf)$ by the following formula:
\chEq\begin{equation}\label{equ:Dehn_is_transv}
(\Dtw{\tmpcurve})_{\hommap}(x) = \intform{[\tmpcurve]}{x}\cdot [\tmpcurve] + x, \qquad \forall x\in H_1(\surf),
\end{equation}
thus it is a transvection along $[\tmpcurve]$, see Eq.~\eqref{equ:transv_def}.

Hence the products of transvections $\aastddif{i}{j}$, $\bbstddif{i}{j}$, $\abstddif{i}{j}$ 
defined by Formula~\eqref{equ:mu_nu_eta} can be realized by products of Dehn twists.
It follows from Theorem~\ref{th:admis_a_b} that all these diffeomorphisms except for $\bbstddif{1}{i}=\bbstddif{i}{1}$ and $\abstddif{i}{1}$ are $\func$-admissible.

On the other hand, by Proposition~\ref{pr:TransAnn_a_1}, $\hdif_{\hommap}$ is generated by the linear isomorphisms $\Dtw{\acurve_i}$, $\Dtw{\bcurve_i}$, $\aastddif{i}{j}$, $\bbstddif{i}{j}$, $\abstddif{i}{j}$, except for $\Dtw{\bcurve_1}$ $\bbstddif{1}{i}=\bbstddif{i}{1}$ and $\abstddif{i}{1}$, where $i\not=j=1,\ldots,\genus.$

Hence, there exists an $\func$-admissible diffeomorphism $\cldif$ of $\surf$
which induces the same isomorphism of $\Hab$ as $\hdif_{\hommap}$.
Then $\trldif = \cldif^{-1}\circ\hdif$ belongs to $\TrlGrM$.
\endproof

\section*{Appendix. Proof of \MainTheorem. Case $\onemanif=\Circle$}
We extend here our proof of \MainTheorem\ given in~\cite{Maks:PhD} to the case when $\manif$ is arbitrary and $\onemanif=\Circle$.

Let $\func,\gfunc:\surf\to\Circle$ be two Morse mappings of same critical type, $\cval$ be their common regular value, $\acurve=\func^{-1}(\cval)$, and $\tmpcurve=\gfunc^{-1}(\cval)$.
By Lemma~\ref{lm:reduce_to_H1_onto} we can assume that homomorphism 
$\func_{\hommap}=\gfunc_{\hommap}:H_1(\surf)\to H_1(\Circle)$ is onto
and by Lemma~\ref{lm:throw_out_1_0_comp} that $\acurve$ and $\tmpcurve$ are connected, i.e.\! \SCCs.

Let us cut $\surf$ along $\acurve$ and denote the obtained surface by $\tsurf$.
Let also $p:\tsurf\to\surf$ be the factor-mapping, $\tfunc:\tsurf\to[0,1]$ the corresponding Morse function induced by $\func$,
$\btsurf{0}=\tfunc^{-1}(0)$, $\btsurf{1}=\tfunc^{-1}(1)$, and $\btsurf{}=\btsurf{0}\cup\btsurf{1}$ (we use the notations of Section~\ref{subsect:CuttingAlongLevelSet}).

\begin{claim}\label{clm:suff_cond_MainTheorem}
If $\acurve=\tmpcurve$, then $\func\sgmh\gfunc$.
\end{claim}
\proof
Since $\func$ and $\gfunc$ are homotopic, we can assume (by small \sgmhom-homotopy) that they coincide near $\acurve$.
Then $\gfunc$ also yields a Morse function $\tgfunc:\surf\to[0,1]$
which coincides near $\btsurf{}$ with $\tfunc$ and $\Crtype{\tfunc}=\Crtype{\tgfunc}$.
By the $\RRR$-case of \MainTheorem\ $\tfunc\sgmh\tgfunc$ with respect to a neighborhood of $\btsurf{}$. Then this \sgmhom-homotopy yields a \sgmhom-homotopy between $\func$ and $\gfunc$ with respect to a neighborhood of $\acurve$.
\endproof

Suppose that $\acurve\not=\tmpcurve$.
Since $\func$ and $\gfunc$ are homotopic, it follows that the restriction $\func|_{\tmpcurve}$ is null-homotopic.
Then by Lemma~\ref{lm:curv_min_int} we can additionally assume that $\acurve\cap\tmpcurve=\varnothing$.

In this case $\ttmpcurve=p^{-1}(\tmpcurve)$ separates $\tsurf$ between $\btsurf{0}$ and $\btsurf{1}$.
Using the method of Claim~\ref{clm:res_min_intersect} we can construct a Morse function $\tfunc_1:\tsurf\to[0,1]$ which coincides with $\tfunc$ near $\btsurf{0}\cup\btsurf{1}$, has critical type of $\tfunc$, and such that $\func_1^{-1}(\frac{1}{2})=\ttmpcurve$.
Then $\tfunc_1$ yields a Morse mapping $\func_1:\surf\to\Circle$ which coincides with $\tfunc$ in a neighborhood of $\acurve$ and such that $\func_1^{-1}(p(\frac{1}{2}))=\tmpcurve$.
Thus $\acurve$ and $\tmpcurve$ are level-sets of $\func_1$.  
Then by Claim~\ref{clm:suff_cond_MainTheorem} we get $\func\sgmh\func_1\sgmh\gfunc$.
\qed


\section{Acknowledgments}
I am sincerely grateful to V.~V.~Sharko, A.~Prishlyak,
M.~Pankov, E.~Po\-lu\-lyakh, I.~Vlasenko for many valuable discussions.
I am deeply thankful to E.~Kudryavtseva for the valuable discussions.
I wish to thank the referee for pointing out to the non-orientable case of \MainTheorem, referring me to the paper~\cite{Kronrod}, and proposing the right name for Kronrod-Reeb graphs.   
I thank A.~Pankov for the information about paper~\cite{Korkmaz}.
I also thank B.~Szepietowski for sending me his paper~\cite{Szepietowski}.


\end{document}